

\def\LaTeX{%
  \let\Begin\begin
  \let\End\end
  \let\salta\relax
  \let\finqui\relax
  \let\futuro\relax}

\def\UK{\def\our{our}\let\sz s}
\def\USA{\def\our{or}\let\sz z}

\UK 



\LaTeX

\USA


\salta

\documentclass[twoside,12pt]{article}
\setlength{\textheight}{24cm}
\setlength{\textwidth}{16cm}
\setlength{\oddsidemargin}{2mm}
\setlength{\evensidemargin}{2mm}
\setlength{\topmargin}{-15mm}
\parskip2mm


\usepackage[usenames,dvipsnames]{color}
\usepackage{amsmath}
\usepackage{amsthm}
\usepackage{amssymb,bbm}
\usepackage[mathcal]{euscript}

\usepackage{cite}
\usepackage{hyperref}
\usepackage{enumitem}

\usepackage[ulem=normalem,draft]{changes}
%
%

%
 
\definecolor{rosso}{rgb}{0.85,0,0}





\bibliographystyle{plain}


%
\newtheorem{theorem}{Theorem}[section]
\newtheorem{remark}[theorem]{Remark}
\newtheorem{corollary}[theorem]{Corollary}

\newtheorem{lemma}[theorem]{Lemma}

\finqui

\def\Beq{\Begin{equation}}

\def\Bthm{\Begin{theorem}}
\def\Ethm{\End{theorem}}
\def\Blem{\Begin{lemma}}
\def\Elem{\End{lemma}}

\def\Brem{\Begin{remark}\rm}
\def\Erem{\End{remark}}

\def\Bdim{\Begin{proof}}
\def\Edim{\End{proof}}
\def\Bcenter{\Begin{center}}
\def\Ecenter{\End{center}}
\let\non\nonumber




\def\step #1 \par{\medskip\noindent{\bf #1.}\quad}
\def\jstep #1: \par {\vspace{2mm}\noindent\underline{\sc #1 :}\par\nobreak\vspace{1mm}\noindent}

\def\Lip{Lip\-schitz}

\def\rhs{right-hand side}




\def\multibold #1{\def\arg{#1}%
  \ifx\arg\pto \let\next\relax
  \else
  \def\next{\expandafter
    \def\csname #1#1#1\endcsname{{\bf #1}}%
    \multibold}%
  \fi \next}

\def\pto{.}

\def\multical #1{\def\arg{#1}%
  \ifx\arg\pto \let\next\relax
  \else
  \def\next{\expandafter
    \def\csname cal#1\endcsname{{\cal #1}}%
    \multical}%
  \fi \next}


\def\multimathop #1 {\def\arg{#1}%
  \ifx\arg\pto \let\next\relax
  \else
  \def\next{\expandafter
    \def\csname #1\endcsname{\mathop{\rm #1}\nolimits}%
    \multimathop}%
  \fi \next}

\multibold
qwertyuiopasdfghjklzxcvbnmQWERTYUIOPASDFGHJKLZXCVBNM.

\multical
QWERTYUIOPASDFGHJKLZXCVBNM.

\multimathop
diag dist div dom mean meas sign supp .


\def\Accorpa #1#2 #3 {\gdef #1{\eqref{#2}--\eqref{#3}}%
  \wlog{}\wlog{\string #1 -> #2 - #3}\wlog{}}


\def\<#1>{\mathopen\langle #1\mathclose\rangle}
\def\norma #1{\mathopen \| #1\mathclose \|}

\def\I2 #1{\int_{Q_t}|{#1}|^2}
\def\IT2 #1{\int_{Q^t}|{#1}|^2}
\def\IO2 #1{\norma{{#1(t)}}^2}

\def\next{\\ & \quad}

\def\iO{\int_\Omega}

\def\dt{\partial_t}
\def\dn{\partial_{\bf n}}
\def\S{{\cal S}}

\def\X{{\cal X}}
\def\Y{{\cal Y}}
\def\Uh{{\cal U}}

\def\checkmmode #1{\relax\ifmmode\hbox{#1}\else{#1}\fi}


\def\bu{{\bf u}}
\def\bv{{\bf v}}
\def\bh{{\bf h}}
\def\bk{{\bf k}}
\def\xih{{\xi^{\bf h}}}
\def\xik{{\xi^{\bf k}}}

\def\bl{{\boldsymbol \lambda}}

\def\erre{{\mathbb{R}}}

\def\enne{{\mathbb{N}}}




\def\genspazio #1#2#3#4#5{#1^{#2}(#5,#4;#3)}
\def\spazio #1#2#3{\genspazio {#1}{#2}{#3}T0}

\def\L {\spazio L}
\def\H {\spazio H}

\def\C #1#2{C^{#1}([0,T];#2)}



\def\Lx #1{L^{#1}(\Omega)}
\def\Hx #1{H^{#1}(\Omega)}

\def\CS #1{C^{#1}(\Sigma)}

\def\Ldue{\Lx 2}

\def\Huno{\Hx 1}
\def\Hdue{\Hx 2}

\def\Liq{{L^\infty(Q)}}




\def\bus{{\bf u}^*}
\let\vp\varphi

\def\a{\alpha}	
\def\b{\beta}	
\def\th{\theta} 
\def\s{\sigma}  
\def\m{\mu}	    
\def\ph{\varphi}	

\def\h{\mathbbm{h}}

\let\TeXchi\chi                         
\newbox\chibox
\setbox0 \hbox{\mathsurround0pt $\TeXchi$}
\setbox\chibox \hbox{\raise\dp0 \box 0 }
\def\chi{\copy\chibox}



\def\uebar{u_1^*}
\def\uzbar{u_2^*}

\def\CP0{(${\mathcal{CP}}_0$)}

\def\CS{{\cal S}}

\def\xih{\xi^{\bf h}}
\def\xik{\xi^{\bf k}}

\let\hat\widehat

\def\uad{{\cal U}_{\rm ad}}
\def\UR{{\cal U}_{R}}
\def\bmu{\mu^*}
\def\bvp{\varphi^*}
\def\bsigma{\sigma^*}

\def\bph{{\ph}^*}
\def\bm{{\m^*}}   
\def\bs{{\s^*}}

\usepackage{amsmath}
\DeclareFontFamily{U}{mathc}{}
\DeclareFontShape{U}{mathc}{m}{it}%
{<->s*[1.03] mathc10}{}

\DeclareMathAlphabet{\mathscr}{U}{mathc}{m}{it}
\Begin{document}


%
\title{Second-order sufficient conditions in the sparse optimal control of a phase field 
    tumor growth \\model with logarithmic potential}
\author{}
\date{}
\maketitle
\Bcenter
\vskip-1.5cm
{\large\sc J\"urgen Sprekels$^{(1)}$}\\
{\normalsize e-mail: {\tt juergen.sprekels@wias-berlin.de}}\\[0.25cm]
{\large\sc Fredi Tr\"oltzsch$^{(2)}$}\\
{\normalsize e-mail:{\tt troeltzsch@math.tu-berlin.de}}\\[0.5cm]
$^{(1)}$
{\small Weierstrass Institute for Applied Analysis and Stochastics}\\
{\small Mohrenstrasse 39, D-10117 Berlin, Germany}\\[0.3cm]
$^{(2)}$
{\small Institut für Mathematik der Technischen Universit\"at Berlin}\\
{\small Strasse des 17. Juni 136, D-10587 Berlin, Germany}\\[10mm]
\Ecenter
\Begin{abstract}
\noindent 
This paper treats a distributed optimal control problem for a tumor growth model of viscous Cahn--Hilliard type. The evolution 
of the tumor fraction is governed by a thermodynamic force induced by a double-well potential of logarithmic type.
The cost functional contains a nondifferentiable
term like the $L^1$--norm in order to enhance the occurrence of sparsity effects in the optimal controls, i.e., of
subdomains of the space-time cylinder where the controls vanish. 
In the context of cancer therapies, sparsity is very important in order that the patient is not exposed to 
unnecessary intensive medical treatment. In this work, we focus on the derivation of second-order sufficient optimality
conditions for the optimal control problem. While in previous works on the system under investigation such conditions 
have been established for the case without sparsity, the case with sparsity has not been treated before.   
\vskip3mm
\noindent {\bf Key words:}
Optimal control, tumor growth models, logarithmic potentials, second-order sufficient optimality conditions, sparsity

\vskip3mm
\noindent {\bf MSC2020 Subject Classification:} {
		35K20, 35K57, 37N25, 49J20, 49J50, 49J52, 49K20, 49K40.		}
\End{abstract}
\salta
\pagestyle{myheadings}
\newcommand\testopari{\sc Sprekels -- Tr\"oltzsch}
\newcommand\testodispari{\sc Second-order conditions in the sparse control of tumor growth}
\markboth{\testopari}{\testodispari}
\finqui
%

\section{Introduction}
\label{INTRO}
\setcounter{equation}{0}

Let  $\a>0,~\b>0,~\chi>0$, and let $\Omega\subset\erre^3$ denote some open and 
bounded domain having a smooth boundary $\Gamma=\partial\Omega$ and the unit outward normal $\,{\bf n}$ with associated
outward normal derivative $\dn$.  Moreover, we fix some final time $T>0$ and
introduce for every $t\in (0,T)$ the sets \,$Q_t:=\Omega\times (0,t)$\, and $\,Q^t:=\Omega\times (t,T)$.
We also set, for convenience,  $Q:=Q_T\,$ and $\,\Sigma:=\Gamma\times (0,T)$.
We then consider the following optimal control problem: 

\vspace{3mm}\noindent
{\bf (CP)} \quad Minimize the cost functional
\begin{align} 
		J((\mu,\vp,\sigma),{\bf u})
   &:=\,  \frac{b_1}2 \iint_Q |\vp-\widehat \vp_Q|^2
  + \frac{b_2}2 \iO |\vp(T)-\widehat\vp_\Omega|^2
  + \frac{b_3}2 \iint_Q (|u_1|^2+|u_2|^2) \nonumber\\
  &\qquad +\,\kappa\,\iint_Q(|u_1|+|u_2|)\nonumber\\
     & =:\, J_1((\mu,\vp,\sigma),{\bf u}) + \kappa g(\bu)
	\label{cost}
\end{align} 
subject to  the state system 
\begin{align}
\label{ss1}
&\alpha\dt\mu+\dt\ph-\Delta\mu=P(\ph)(\sigma+\chi(1-\ph)-\mu) - \h(x,t)
u_1 \quad&&\mbox{in }\,Q\,,\\[1mm]
\label{ss2}
&\beta\dt\vp-\Delta\vp+F_1'(\vp) + F_2'(\vp)=\mu+\chi\,\sigma \quad
&&\mbox{in }\,Q\,,\\[1mm]
\label{ss3}
&\dt\sigma-\Delta\sigma=-\chi\Delta\vp-P(\ph)(\sigma+\chi(1-\ph)-\mu)+u_2\quad&&\mbox{in }\,Q\,,\\[1mm]
\label{ss4}
&\dn \mu=\dn\vp=\dn\sigma=0 \quad&&\mbox{on }\,\Sigma\,,\\[1mm] 
\label{ss5}
&\mu(0)=\mu_0,\quad \vp(0)=\vp_0,\quad \sigma(0)=\sigma_0\,\quad &&\mbox{in }\,\Omega\,,
\end{align}
\Accorpa\Statesys {ss1} {ss5}
and to the control constraint
\begin{equation}
\label{cont:constr}
{\bf u}=(u_1,u_2)\in\uad\,.
\end{equation}
Here, $b_1\,$ and $\,b_2\,$ are nonnegative constants, while $\,b_3\,$ and $\,\kappa\,$  
are positive;  $\widehat\vp_Q$ and $\widehat \vp_\Omega$ are given target functions. The term $g(\bu)$ accounts for possible sparsity effects.
 Moreover, the set of admissible controls $\uad$ is a nonempty, closed and convex subset of the control space
\begin{equation}
\label{defU}
\Uh:= L^\infty(Q)^2.
\end{equation}

The state system \Statesys\ constitutes a simplified and relaxed version of the four-species thermodynamically consistent model for tumor growth originally proposed by Hawkins-Daruud 
et al.\ in \cite{HZO} that additionally includes the chemotaxis-like terms $\,\chi\sigma\,$ in \eqref{ss2} and 
$\,-\chi\Delta\vp\,$ in \eqref{ss3}. 
Let us briefly review the role of the occurring symbols. The primary (state) variables $\ph, \m, \s$ 
denote the tumor fraction, the associated chemical potential, and the nutrient concentration, respectively.
Furthermore, the additional term $\a\dt\m$ corresponds to a parabolic regularization of equation \eqref{ss1},
while $\b\dt\ph$ is the viscosity contribution to the Cahn--Hilliard equation.
The nonlinearity $P$ denotes a proliferation function, whereas the positive constant $\chi$
represents the chemotactic sensitivity
and provides the system with a cross-diffusion coupling.

The evolution of the tumor fraction is mainly governed by
the nonlinearities $F_1$ and $F_2$ whose derivatives occur in \eqref{ss2}. Here, $F_2$ is smooth, typically a 
concave function. As far as
$F_1$ is concerned, we admit in this paper functions of logarithmic type such as   
\begin{align}
\label{logpot} 
&F_{\rm 1,log}(r)=\left\{\begin{array}{ll}
(1+r)\,\ln(1+r)+(1-r)\,\ln(1-r)\quad&\mbox{for $\,r\in (-1,1)$}\\
2\,\ln(2)\quad&\mbox{for $\,r\in\{-1,1\}$\,,}\\
+\infty\quad&\mbox{for $\,r\not\in[-1,1]$}
\end{array}\right.
\end{align}
We assume that $F=F_1+F_2$ is a double-well potential. This is actually the case 
if $F_2(r)=k(1-r^2)$ with a sufficiently large $\,k>0$. Note also that $F'_{\rm 1,log}(r)$ becomes 
unbounded as $r\searrow -1$ and $r\nearrow 1$.

The control variable $u_2$ occurring in \eqref{ss3} can model  
either an external medication or some nutrient supply, while
$u_1$, which occurs in the phase equation \eqref{ss1}, models the application of a cytotoxic drug to 
the system. Usually, $u_1\,$ is multiplied by a truncation function $\,\h(\varphi)\,$ in order to have the action
only in the spatial region where the tumor cells are located. Typically,
one assumes that $\h(-1)=0, ~\h(1)=1$, and $\h(\varphi)$ is in between if $-1<\ph<1$;
see \cite{GLSS, GARL_1, HKNZ, KL} for some insights on possible choices of $\h$. Also in \cite{CSS1,CSS2,ST},
this kind of nonlinear coupling between $\,u_1\,$ and $\,\varphi\,$ has been admitted. For our following analysis,
this nonlinear coupling is too strong, and, only for technical reasons, we have chosen to simplify the original system somewhat by 
assuming that $\,\h=\h(x,t)\,$ is a bounded nonnegative function that does not depend on $\,\varphi$. 
We stress the fact that this simplification does not 
have any impact on the validity of the results from \cite{CSS2} to be used below.

As far as well-posedness {is concerned}, the above model was already investigated in the 
case $\chi=0$ in \cite{CGH,CGRS1,CGRS2,CRW},
and in \cite{FGR} with $\a=\b=\chi=0$.
There the authors also pointed out how the relaxation parameters $\a$ and $\b$ can be set to zero, by providing
the proper framework in which a limit system can be identified and uniquely solved.
We also note that in \cite{CGS24} a version has been studied in which the Laplacian in the equations 
\eqref{ss1}--\eqref{ss3} has been replaced by fractional powers of a more general class of selfadjoint operators
having compact resolvents. A model which is similar to the one studied in this note was the subject
of \cite{CSS1,ST}.

For some nonlocal variations of the above model we refer to \cite{FLRS, FLS, SS}.
Moreover, in order to better emulate in-vivo tumor growth,
it is possible to include in similar models the effects generated by the fluid flow development
by postulating a Darcy law or a Stokes--Brinkman law.
In this direction, we refer to \cite{WLFC,GLSS,DFRGM ,GARL_1,GARL_4,GAR, EGAR, FLRS,GARL_2, GARL_3, CGSS1}, and {we also mention} \cite{GLS}, where elastic effects are included.
For further models, discussing the case of multispecies,
we refer the reader to \cite{DFRGM,FLRS}. 
The investigation of associated optimal control problems also presents a 
wide number of results of which we mention{\cite{CGRS3,EK, EK_ADV, S, S_a, S_b, S_DQ, SigTime, FLS, ST, CGS24, CSS1, GARLR, 
KL, SW, GLS_OPT}. 

Sparsity in the optimal control theory of partial differential equations is a very active field of research.  
The use of sparsity-enhancing functionals goes back to inverse problems and image processing. Soon after
the seminal paper  \cite{stadler2009}, many results were published. 
We mention only very few of them with closer relation to our paper, in particular 
\cite{casas_herzog_wachsmuth2017,herzog_obermeier_wachsmuth2015,herzog_stadler_wachsmuth2012}, on directional sparsity,
 and \cite{casas_troeltzsch2012} on a general theorem for second-order conditions;
moreover, we refer to some new trends in the investigation of sparsity, namely, infinite horizon sparse 
optimal control (see, e.g., \cite{Kalise_Kunisch_Rao2017,Kalise_Kunisch_Rao2020}), and fractional order optimal 
control  (cf. \cite{Otarola_Salgado2018}, \cite{Otarola2020}). 
These papers concentrated on first-order optimality conditions for sparse optimal controls of single elliptic and parabolic equations. 
In  \cite{casas_ryll_troeltzsch2013,casas_ryll_troeltzsch2015}, first- and second-order optimality conditions 
have been discussed in the context of sparsity for the (semilinear) system of  FitzHugh--Nagumo equations. 
Moreover, we refer to the measure control of the Navier--Stokes system 
studied in \cite{Casas_Kunisch2021}.

The optimal control problem {\bf (CP)}  has recently been investigated in \cite{CSS2} for the case of 
logarithmic potentials $F_1$ and without sparsity terms, where second-order sufficient optimality conditions have 
been derived using the $\,\tau$--critical cone and the splitting technique as described in the textbook \cite{Fredibuch}. 
In \cite{ST} and \cite{CSS4}, sparsity terms have been incorporated, where in the latter paper not only 
logarithmic nonlinearities but also nondifferentiable double obstacle potentials have been admitted. However, second-order sufficient optimality conditions have not been derived.

The derivation of meaningful second-order conditions for locally  optimal controls of {\bf (CP)}
in the logarithmic case with sparsity term is the main object of this paper. In particular, we aim at constructing suitable
critical cones which are as small as possible. In our approach, we follow the recent work \cite{SpTr2} on the sparse
optimal control of Allen--Cahn systems, which was based on ideas developed in \cite{casas_ryll_troeltzsch2015}.  

The paper is organized as follows. 
In the next section, we list and discuss our assumptions, and we collect known results from \cite{CSS4}	concerning the properties of the
state system \Statesys\ and of the control-to-state operator.
In Section 3, we study the optimal control problem. We derive first-order necessary optimality conditions 
and results concerning the full sparsity of local minimizers, and we establish second-order sufficient optimality conditions for the optimal control problem {\bf (CP)}. 
In an appendix, we prove auxiliary results that are needed for the main theorem on second-order sufficient conditions.

Prior to this, let us fix some notation.
For any Banach space $\,X$, we denote by \linebreak$\|\cdot\|_X$
the norm in the space $\,X$, by $\,X^*\,$ its dual space, 
and by $\langle \, \cdot\, , \, \cdot \, \rangle_{X}$ the duality pairing between $\,X^*\,$ and $\,X$. 
For any $1 \leq p \leq \infty$ and $k \geq 0$, we denote the standard Lebesgue and Sobolev spaces on $\,\Omega\,$ 
by $L^p(\Omega)$ and $W^{k,p}(\Omega)$, and the corresponding norms by 
$\norma{\,\cdot\,}_{L^p(\Omega)}=\norma{\,\cdot\,}_{p}$ and $\norma{\,\cdot\,}_{W^{k,p}(\Omega)}$, respectively. 
For $p = 2$, they become Hilbert spaces, and we employ the standard notation $H^k(\Omega) := W^{k,2}(\Omega)$. 
As usual, for Banach spaces $\,X\,$ and $\,Y\,$ that are both continuously embedded in some topological vector space~$Z$,
we introduce the linear space
$\,X\cap Y\,$ which becomes a Banach space when equipped with its natural norm $\,\|v\|_{X \cap Y}:=
\|v\|_X\,+\,\|v\|_Y$, for $\,v\in X\cap Y$.
Moreover, we recall the definition~\eqref{defU} of the control space ${\cal U}$ and introduce the spaces
\begin{align}
  & H := \Ldue \,, \quad  
  V := \Huno\,,   \quad
  W_{0} := \{v\in\Hdue: \ \dn v=0 \,\mbox{ on $\,\Gamma$}\}.
  \label{defHVW}
\end{align}
Furthermore, by $\,(\,\cdot\,,\,\cdot\,)$ and  $\,\Vert\,\cdot\,\Vert$ we denote the standard inner product 
and related norm in $\,H$, and, for simplicity,  we also set $\langle\,\cdot\,,\,\cdot\,\rangle := 
\langle\,\cdot\,,\,\cdot\,\rangle_V $.

Throughout the paper, we make repeated use of H\"older's inequality, of the elementary Young inequality
\begin{equation}
\label{Young}
a b\,\le \delta |a|^2+\frac 1{4\delta}|b|^2\quad\forall\,a,b\in\erre, \quad\forall\,\delta>0,
\end{equation}
as well as the continuity of the embeddings $H^1(\Omega)\subset L^p(\Omega)$ for $1\le p\le 6$ and 
$\Hdue\subset C^0(\overline\Omega)$.

We close this section by introducing a convention concerning the constants used in estimates within this paper: we denote by $\,C\,$ any 
positive constant that depends only on the given data occurring in the state system and in the cost functional, as well as 
on a constant that bounds the $\left(L^\infty(Q)\times L^\infty(Q)\right)$--norms of the elements of $\uad$. The actual value of 
such generic constants $\,C\,$ 	may 
change from formula to formula or even within formulas. Finally, the notation $C_\delta$ indicates a positive constant that
additionally depends on the quantity $\delta$.


\section{General setting  and properties of the control-to-state operator}
\label{STATE}
\setcounter{equation}{0}

In this section, we introduce the general setting of our control 
problem and state some results on the state system \eqref{ss1}--\eqref{ss5} and the control-to-state operator
 that in its present form have been
established in \cite{CSS2,CSS4}.

We make the following assumptions on the data of the system.

\begin{description}
\item[(A1)] \quad $\alpha,\,\beta,\,\chi\,$ are positive constants.
\item[(A2)] \quad $F=F_1+F_2$, where $F_2 \in C^5(\erre)$ has a \Lip\ continuous derivative $F'_2$, and where $\,F_1:\erre\to [0,+\infty]\,$ is convex and lower semicontinuous and satisfies $\,F_1(0)=0$,  $\,F_{1_{|(-1,1)}}\in C^5(-1,1)$, as well as 
\begin{equation}
\label{singular}
\lim_{r\searrow -1} F_1'(r)=-\infty \quad\mbox{ and }\quad \lim_{r\nearrow 1} F_1'(r)=+\infty\,.
\end{equation}
\item[(A3)] \quad $P \in C^3(\erre)\cap W^{3,\infty}(\erre)\,$ and $\,\h \in \Liq\,$  are nonnegative and bounded.
\item[(A4)] \quad The initial data satisfy $\,\mu_0,\sigma_0\in \Huno\cap L^\infty(\Omega)$, $\vp_0\in W_0$, as well as
\begin{equation}
\label{sepini}
-1 < \min_{x\in \overline\Omega}\,\vp_0(x) \le \max_{x\in\overline\Omega}\,\vp_0(x) < 1\,.
\end{equation}
\item[(A5)] \quad 	With fixed given constants $\underline u_i,\overline u_i$ satisfying $\underline u_i
<\overline u_i$, $i=1,2$, we have	
\begin{equation}
\label{defuad}
\uad=\left\{\bu=(u_1,u_2)\in {\cal U}: \underline u_i\le u_i\le\overline u_i\,\mbox{ a.e. in }\,Q\,\mbox{ for }\,i=1,2\right\}. 
\end{equation}
\item[(A6)] \quad $R>0$ is a constant such that $\,\uad\subset{\cal U}_R:=\{\bu\in{\cal U}:\,\|\bu\|_{\cal U}<R\}$.
\end{description}

\Brem		
Observe that {\bf (A3)} implies that the functions $P,P',P''$ are \Lip\ continuous on $\erre$.
Let us also note that  $F_1=F_{\rm 1,log}$ satisfies {\bf (A2)}. Moreover, \eqref{sepini} implies that 
initially there are no pure phases. Finally, {\bf (A6)} just fixes an open and bounded subset of $\,{\cal U}\,$
that contains $\,\uad$. 
\Erem 

The following result is a consequence of \cite[Thm.~2.3]{CSS4}.
\Bthm
\label{THM:STRONG}
Suppose that the conditions {\bf (A1)}--{\bf (A6)} are fulfilled. 
Then the state system \eqref{ss1}--\eqref{ss5} has for every $\bu=(u_1,u_2)\in \UR$ a unique strong solution
$(\mu,\vp,\sigma)$ with the regularity
\begin{align}
\label{regmu}
&\mu\in  H^1(0,T;H) \cap C^0([0,T];V) \cap L^2(0,T;W_0)\cap  L^\infty(Q),\\[1mm]
\label{regphi}
&\ph\in W^{1,\infty}(0,T;H)\cap H^1(0,T;V)\cap L^\infty(0,T;W_0) \cap C^0(\overline Q),
\\[1mm]
\label{regsigma}
&\sigma\in H^1(0,T;H)\cap C^0([0,T];V)\cap L^2(0,T;W_0)\cap L^\infty(Q).
\end{align}
Moreover, there is a constant $K_1>0$, which depends on $\Omega,T,R,\alpha,\beta$ and the data of the 
system, but not on the choice of $\bu\in \UR $, such that
\begin{align}\label{ssbound1}
&\|\mu\|_{H^1(0,T;H) \cap C^0([0,T];V) \cap L^2(0,T;W_0)\cap L^\infty(Q)}\nonumber\\[1mm]
&+\,\|\ph\|_{W^{1,\infty}(0,T;H)\cap H^1(0,T;V) \cap L^\infty(0,T;W_0)\cap C^0(\overline Q)}
\nonumber\\[1mm]&+\,\|\sigma\|_{H^1(0,T;H) \cap C^0([0,T];V) \cap L^2(0,T;W_0)\cap L^\infty(Q)}\,\le\,K_1\,.
\end{align}
Furthermore, there are constants $r_*,r^*$, which depend on $\Omega,T,R,\alpha,\beta$ and the data of the 
system, but not on the choice of $\bu\in \UR$, such that
\begin{equation}\label{ssbound2}
-1  <r_*\le\vp(x,t)\le r^*<1 \quad\mbox{for all $(x,t)\in \overline Q$}.
\end{equation}
Also, there is some constant $K_2>0$ having the same dependencies as $K_1$ such that
\begin{align}
\label{ssbound3}
	&\max_{i=0,1,2,3}\,\left\|P^{(i)}(\vp)\right\|_{L^\infty(Q)}\,
	+\,\max_{i=0,1,2,3,4{,5}}\, \max_{j=1,2}\,\left\|
	F_j^{(i)}(\vp)\right\|_{L^\infty(Q)} \,\le\,K_2\,.
\end{align}
Finally, if $\,\bu_i\in\UR\,$ are given controls and $(\mu_i,\vp_i,\sigma_i)$ the corresponding solutions
to \Statesys, for $\,i=1,2$, then, with a constant $\,K_3>0\,$ having the same dependencies as $\,K_1$,
\begin{align}
\label{stabu}
&\|\mu_1-\mu_2\|_{H^1(0,T;H)\cap C^0([0,T];V)\cap L^2(0,T;W_0)}\,+\,
\|\vp_1-\vp_2\|_{H^1(0,T;H)\cap C^0([0,T];V)\cap L^2(0,T;W_0)}\nonumber\\
&+\,\|\s_1-\s_2\|_{H^1(0,T;H)\cap C^0([0,T];V)\cap L^2(0,T;W_0)}\,\le\,K_3\,\|\bu_1-\bu_2\|_{L^2(Q)^2}\,.
\end{align} 
\Ethm

\Brem
Condition \eqref{ssbound2}, known as the {\em separation property},
is especially important for the case of singular 
potentials such as $F_1=F_{\rm 1,log}$, since it guarantees that the phase variable \,$\ph$\, always stays away
from the critical values $\,-1,+1$. The singularity of $\,F_1'\,$ is therefore no longer an obstacle for the analysis,
as the values of $\ph$ range in some interval in which 
$\,F_1'\,$ is smooth.
\Erem         

Owing to Theorem 2.2, the control-to-state operator $$\, {\cal S}:\bu=(u_1,u_2)\mapsto 
(\mu,\ph,\sigma)\,$$ is well defined as a
mapping between ${\cal U}=L^\infty(Q)^2$ and the Banach space specified by the regularity 
results~\eqref{regmu}--\eqref{regsigma}. We now discuss its differentiability properties. 
For this purpose, some functional analytic preparations are in order. We first define the linear spaces
\begin{align}
\X\,&:=\,X\times \widetilde X\times X, \mbox{\,\,\, where }\nonumber\\
X\,&:=\,H^1(0,T;H)\cap C^0([0,T];V)\cap L^2(0,T;W_0)\cap L^\infty(Q), \nonumber\\
\widetilde X\,&:=W^{1,\infty}(0,T;H)\cap H^1(0,T;V)\cap L^\infty(0,T;W_0)\cap C^0(\overline Q),
\label{calX}
\end{align}
which are Banach spaces when endowed with their natural norms. Next, we introduce the linear space
\begin{align}
&\Y\,:=\,\bigl\{(\mu,\ph,\sigma)\in \calX: \,\alpha\dt\mu+\dt\ph-\Delta\mu\in \Liq, 
\,\,\,\beta\dt\ph-\Delta\ph-\mu\in \Liq,\nonumber\\
& \hspace*{14mm}\dt\sigma-\Delta\sigma+\chi\Delta\ph\in\Liq\bigr\},
\label{calY}
\end{align}
which becomes a Banach space when endowed with the norm
\begin{align}
\|(\mu,\ph,\sigma)\|_{\Y}\,:=\,&\|(\mu,\ph,\sigma)\|_\calX
\,+\,\|\alpha\dt\mu+\dt\ph-\Delta\mu\|_{\Liq}
\,+\,\|\beta\dt\ph-\Delta\ph-\mu\|_{\Liq}\nonumber\\
&+\,\|\dt\sigma-\Delta\sigma+\chi\Delta\ph\|_{\Liq}\,.
\label{normY}
\end{align}
Finally, we put 
\begin{align}
&Z:=H^1(0,T;H)\cap C^0([0,T];V)\cap L^2(0,T;W_0),
\label{defZ}
\\
&{\cal Z}:= Z\times Z\times Z.
\label{calZ}
\end{align} 

For fixed $\,(\vp^*,\mu^*,\s^*)$, we first discuss an auxiliary result for the linear initial-boundary value problem
\begin{align}\label{aux1}
&\alpha\dt\mu+\dt\varphi-\Delta\mu\,=\,\lambda_1\left[P(\bvp)(\sigma-\chi\varphi-\mu)+P'(\bvp)(\bsigma+\chi(1-\bvp)-\bmu)\ph
\right]\nonumber\\
&\hspace*{37mm}-\lambda_2\, \h\, h_1+\lambda_3 f_1\quad\mbox{in }\,Q,\\[1mm]
\label{aux2}
&\beta\dt\varphi-\Delta\varphi-\mu\,=\,\lambda_1\left[\chi\,\sigma-F''(\bvp)\varphi\right]+\lambda_3 f_2 \quad\mbox{ in \,$Q$},
\\[1mm]
\label{aux3}
&\dt\sigma-\Delta\sigma+\chi\Delta\varphi\,=\,\lambda_1\left[-P(\bvp)(\sigma-\chi\varphi-\mu)-
P'(\bvp)(\bsigma+\chi(1-\bvp)-\bmu)\varphi \right]\nonumber\\
&\hspace*{37mm}+\lambda_2 h_2+\lambda_3 f_3\quad\mbox{ in \,$Q$},\\[1mm]
\label{aux4}
&\dn\mu\,=\,\dn\varphi\,=\,\dn\sigma\,=\,0 \quad\mbox{ on \,$\Sigma$},\\[1mm]
\label{aux5}
&\mu(0)\,=\,\lambda_4 \mu_0,\quad \varphi(0)\,=\,\lambda_4\varphi_0, \quad \sigma(0)\,=\,\lambda_4
\sigma_0,\,\,\,\mbox{ in }\,\Omega,
\end{align}
which for $\lambda_1 = \lambda_2 = 1$ and $\lambda_3 = \lambda_4 = 0$ coincides  with the linearization 
of the state equation at $((\bmu,\bvp,\bsigma),(\uebar,\uzbar))$. 
We have the following result.
\Blem \label{Hilfe}
Suppose that $\lambda_1,\lambda_2,\lambda_3,\lambda_4 \in\{0,1\}$ are given
and that the assumptions {\bf (A1)}--{\bf (A6)} are fulfilled. Moreover, let $\bus=(u_1^*,u_2^*)\in \UR$ be given and
$ (\mu^*,\bvp,\bsigma)=\CS(\bus)$. Then \eqref{aux1}--\eqref{aux5}
has for every $\bh=(h_1,h_2)\in L^2(Q)^2$ and $(f_1,f_2,f_3)\in L^2(Q)^3$ a unique solution
\,$(\mu,\vp,\sigma)\in Z\times \widetilde X\times Z$. 
Moreover, the linear mapping 
\begin{equation}
\label{linsol}
((h_1,h_2),(f_1,f_2,f_3))\mapsto (\mu,\ph,\sigma)
\end{equation}
is continuous from $\,L^2(Q)^2\times L^2(Q)^3\,$ into $Z\times\widetilde X\times Z$. Moreover, if 
$\bh\in \Liq^2$ and $(f_1,f_2,f_3)\in\Liq^3$, in addition, then it holds
$(\mu,\varphi,\sigma)\in{\cal Y}$, and the mapping \eqref{linsol} is continuous from  $\Liq^2 \times \Liq^3$ into
${\cal Y}$.  
\Elem                 
\Bdim
The existence result and the continuity of the mapping \eqref{linsol} between the spaces $\,\Liq^2\times\Liq^3\,$ and $\,{\cal Y}\,$
directly follow from the statement of \cite[Lem.~4.1~and~Rem.~4.2]{CSS2}. Moreover, from the estimates (4.36)--(4.38) and (4.43) 
in \cite{CSS2} we can conclude that the mapping \eqref{linsol} is also continuous between the spaces
$\,L^2(Q)^2\times L^2(Q)^3\,$ and $\,Z\times\widetilde X\times Z$. 
\Edim
Now let $\bus=(\uebar,\uzbar)\in {\cal U}_R$ be arbitrary and 
$(\bmu,\bvp,\bsigma)={\cal S}(\bus)$. Then, according to \cite[Thm.~4.4]{CSS2}, the control-to-state 
operator $\S$ is twice continuously Fr\'echet differentiable at $\,\bus\,$ as a mapping from $\,{\cal U}\,$ into 
$\,{\cal Y}$. Moreover, for every $\,\bh=(h_1,h_2)\in {\cal U}$, the first Fr\'echet derivative $\,S'(\bus)\in 
{\cal L}({\cal U},\Y)\,$ of $\,\S\,$ at
$\,\bus\,$ is given by the identity $\,\S'(\bus)[\bh]=(\eta^\bh,\xi^\bh,\theta^\bh)$, where 
$\,(\eta^\bh,\xi^\bh,\theta^\bh)\in {\cal Y}\,$ is the 
unique solution to the linearization of the state system given by the initial-boundary value problem \eqref{aux1}--\eqref{aux5}
with $\lambda_1 = \lambda_2 = 1$ and $\lambda_3 = \lambda_4 = 0$. 

\vspace*{2mm}
\Brem
Observe that, in view of the continuity of the embedding $\,{\cal Y}\subset Z\times\widetilde X\times  Z$, the operator 
$\,\S'(\bus)\in {\cal L}(\Uh,{\cal Y})\,$ also belongs to the space $\,{\cal L}(\Uh,Z\times\widetilde X\times  Z)$ and,
owing to the density of $\,\Uh\,$ in $\,L^2(Q)^2$, can be extended continuously to an element
of ${\cal L}(L^2(Q)^2,Z\times\widetilde X\times Z)$ without changing its operator norm. Denoting the extended operator
still by $\,\S'(\bus)$, we see 
that the identity $\S'(\bus)[\bh]=(\eta^\bh,\xih,\theta^\bh)$ is also valid
for every $\,\bh\in L^2(Q)^2$, only that $\,(\eta^\bh,\xih,\theta^\bh)\in Z\times\widetilde X\times  Z$, in general. 
In addition, it also follows from the proof of \cite[Lem.~4.1]{CSS2} that
there is a constant $\,K_4>0$, which depends only on $\,R\,$ and the data, such that
\begin{equation}
\label{lsbound1}
\|\S'(\bu)[\bh]\|_{Z\times\widetilde X\times  Z}\,\le\,K_4\,\|\bh\|_{L^2(Q)^2}\,\quad\mbox{for all 
\,$\bu\in\UR$\, and every $\,\bh\in L^2(Q)^2$}\,.
\end{equation}
\Erem 

\vspace*{2mm}
Next, we show a Lipschitz property for the extended operator $\,\S'$.
\Blem
The mapping $\,\S':{\cal U}\to {\cal L}(L^2(Q)^2,Z\times\widetilde X\times Z)$, $\bu\mapsto \S'(\bu)$, is Lipschitz continuous in the
following sense: there is a constant $K_5>0$, which depends only on $\,R\,$ and the data, such that, for all controls $\,\bu_1,\bu_2
\in\UR\,$ and all increments $\bh\in L^2(Q)^2$,
\begin{align}
\label{Lip1}
\|\left(\S'(\bu_1)-\S'(\bu_2)\right)[\bh]\|_{\cal Z}\,\le\,K_5\,\|\bu_1-\bu_2\|_{L^2(Q)^2}\,\|\bh\|_{L^2(Q)^2}\,.
\end{align}
\Elem
\Bdim
We put $\,(\mu_i,\vp_i,\s_i):=\S(\bu_i)$, $\,(\eta_i,\xi_i,\theta_i):=\S'(\bu_i)[\bh]$, $i=1,2$, as well as
\begin{align*}
&\bu:=\bu_1-\bu_2, \quad \mu:=\mu_1-\mu_2,\quad \vp:=\vp_1-\vp_2,\quad \s:=\s_1-\s_2,\\
&\eta:=\eta_1-\eta_2,\quad\xi:=\xi_1-\xi_2,\quad \theta:=\theta_1-\theta_2\,.
\end{align*} 
Then it follows from \eqref{stabu} in Theorem 2.2 that
\begin{equation}
\label{lippi0}
\|(\mu,\vp,\s)\|_{{\cal Z}}\,\le\,K_3\|\bu\|_{L^2(Q)^2}\,.
\end{equation}
Moreover, $\,(\eta,\xi,\theta)\,$ solves the problem
\begin{align}
\label{lippi1}
&\alpha \dt\eta+\dt\xi-\Delta\eta\,=\,P(\vp_1)(\theta-\chi\xi-\eta)+P'(\vp_1)(\s_1+\chi(1-\vp_1)-\mu_1)\xi\,+\,f_1\,,\\
\label{lippi2}
&\beta\dt\xi-\Delta\xi=\chi\theta-F''(\vp_1)\xi\,+\,f_2\,,\\
\label{lippi3}
&\dt\theta-\Delta\theta+\chi\Delta\xi=-P(\vp_1)(\theta-\chi\xi-\eta)-P'(\vp_1)(\s_1+\chi(1-\vp_1)-\mu_1)\xi\,+\,f_3\,,\\
\label{lippi4}
&\dn\eta=\dn\xi=\dn\theta=0\,,\\
\label{lippi5}
&\eta(0)=\xi(0)=\theta(0)=0\,,
\end{align}
which is of the form \eqref{aux1}--\eqref{aux5} with $\lambda_1=\lambda_3=1$ and $\lambda_2=\lambda_4=0$, and where
\begin{align}
\label{fiete1}
&f_1:=-f_3:=((P(\vp_1)-P(\vp_2))(\theta_2-\chi\xi_2-\eta_2)\,+\,P'(\vp_1)(\s-\chi\vp-\mu)\xi_2\nonumber\\
&\hspace*{23mm}
+(P'(\vp_1)-P'(\vp_2))(\s_2+\chi(1-\vp_2)-\mu_2)\xi_2\,,\\
\label{fiete2}
&f_2:= -(F''(\vp_1)-F''(\vp_2))\xi_2\,.
\end{align}
We therefore conclude from Lemma 2.4 that
\begin{equation*}
\|(\eta,\xi,\theta)\|_{\cal Z}\,\le\,C\left(\|f_1\|_{L^2(Q)}\,+\,\|f_2\|_{L^2(Q)}\right)\,.
\end{equation*}
Hence, the proof will be finished once we can show that
\begin{equation}
\label{fiete3}
\|f_1\|_{L^2(Q)}+\|f_2\|_{L^2(Q)}\,\le\,C\,\|\bu\|_{L^2(Q)^2}\,\|\bh\|_{L^2(Q)^2}\,.
\end{equation}
To this end, we first use the mean value theorem, \eqref{ssbound3}, H\"older's inequality, the continuity of the
embedding $\,V\subset L^4(\Omega)$, as well as \eqref{stabu} and \eqref{lippi0}, to find that
\begin{align}
\label{fiete4}
&\|f_2\|_{L^2(Q)}^2\,\le\,C\!\iint_Q \!|\vp|^2\,|\xi_2|^2\,\le\,C\!\int_0^T \! \|\vp\|_4^2\,\|\xi_2\|_4^2\,ds
\,\le\,C\,\|\vp\|_{C^0([0,T];V)}^2\,\|\xi_2\|_{C^0([0,T];V)}^2
\nonumber\\[1mm]
&\le\,C\,\|\vp\|_Z^2\,\|\S'(\bu_2)[\bh]\|_{{\cal Z}}^2
\,\le\,C\,\|\bu\|_{L^2(Q)^2}^2\,\|\bh\|_{L^2(Q)^2}^2\,.
\end{align}
Here, we have for convenience omitted the argument $\,s\,$ in the third integral. We will do this repeatedly in the
following.
For the three summands on the right-hand side of \eqref{fiete1}, which we denote by $\,A_1, A_2, A_3$, in this order, we obtain by
similar reasoning the estimates
\begin{align}
\label{fiete5}
&\iint_Q|A_1|^2\,\le\,C\iint_Q |\vp|^2\,|\theta_2-\chi\xi_2-\eta_2|^2\,\le\,C\int_0^T \|\vp\|_4^2\left(\|\theta_2\|_4^2
\,+\,\|\xi_2\|^2_4\,+\,\|\eta_2\|_4^2\right)ds \nonumber\\[1mm]
&\quad\le\,C\,\|\bu\|_{L^2(Q)^2}^2\,\|\bh\|_{L^2(Q)^2}^2\,,\\[2mm]
\label{fiete6}
&\iint_Q|A_2|^2\,\le\,C\int_0^T |\xi_2|^2 \left(|\mu|^2+|\vp|^2+|\s|^2\right)ds\,\le\,C\int_0^T\|\xi_2\|_4^2\left(
\|\mu\|_4^2+\|\vp\|_4^2+\|\s\|_4^2\right)ds\nonumber\\[1mm]
&\quad\le C\,\|\bu\|_{L^2(Q)^2}^2\,\|\bh\|_{L^2(Q)^2}^2\,,\\[2mm]
\label{fiete7}
&\iint_Q|A_3|^2\,\le\,C\left(\|\s_2\|_{\Liq}^2+\|\vp_2\|_{\Liq}^2+\|\mu_2\|_{\Liq}^2+1\right)\iint_Q|\vp|^2|\xi_2|^2
\nonumber\\[1mm]
&\quad\le\,C\,\|\bu\|_{L^2(Q)^2}^2\,\|\bh\|_{L^2(Q)^2}^2\,,
\end{align}
where in the last estimate we also used \eqref{ssbound1} and \eqref{fiete4}. With this, the assertion is proved. 
\Edim

\vspace*{2mm}
Next, we turn our interest to the second Fr\'echet derivative $\,\S''(\bus)\,$ of $\,\S\,$ at $\,\bus$. Let
$\,\bh=(h_1,h_2)\in\Uh\,$ and $\,\bk=(k_1,k_2)\in\Uh$. Then, $\,(\eta^\bh,\xih,\theta^\bh):=\S'(\bus)[\bh]\,$
and $\,(\eta^\bk,\xik,\theta^\bk):=\S'(\bus)[\bk]$ both belong to $\,{\cal Y}$ and, by virtue of 
\cite[Thm.~4.6]{CSS2}, $\,(\nu,\psi,\rho)=
\S''(\bus)[\bh,\bk]\in {\cal Y}\,$ is the unique solution to the bilinearization of the state system at 
$((\bmu,\bvp,\bsigma),(\uebar,\uzbar))$, which is given by the 
linear initial-boundary value problem
\begin{align}
\label{bilin1}
&\alpha\dt\nu+\dt\psi-\Delta\nu\,= P(\vp^*)(\rho-\chi\psi-\nu)+P'(\vp^*)(\s^*+\chi(1-\vp^*)-\mu^*)\psi
+ f_1  \,,\\[1mm]
\label{bilin2}
&\beta\dt\psi-\Delta\psi-\nu = \chi \rho  - F''(\bvp)\psi +f_2 \,,
\\[1mm]
\label{bilin3}
&\dt\rho-\Delta\rho+\chi\Delta\psi\,=\,- P(\vp^*)(\rho-\chi\psi-\nu)-P'(\vp^*)(\s^*+\chi(1-\vp^*)-\mu^*)\psi+f_3\,,\\[1mm]
\label{bilin4}
&\dn\nu\,=\,\dn\psi\,=\,\dn\rho\,=\,0\,,\\[1mm]
\label{bilin5}
&\nu(0)\,=\, \psi(0)\,=\, \rho(0)\,=\,0\,,
\end{align}
\Accorpa\Bilin {bilin1} {bilin5}
and which is again of the form \eqref{aux1}--\eqref{aux5} with $\,\lambda_1=\lambda_3=1\,$ and $\,\lambda_2=\lambda_4=0$, 
where 
\begin{align}
\label{fiete8}
	&f_1:=-f_3:= P'(\bph) \left(\xik\,(\theta^\bh-\chi\xih-\eta^\bh)+\xih \,(\theta^\bk - \chi \xik - \eta^\bk)\right)\nonumber\\
	&\hspace*{25mm} + P''(\bph)\,\xik\,\xih\,(\s^* + \chi (1-\vp^*) - \mu^*)\,,\\[2mm]
\label{fiete9}	
&f_2:= -F^{(3)}(\vp^*)	\xih\xik\,.
\end{align}

Now assume that  $\,\bh,\bk\in L^2(Q)^2$ are given. Then the expressions $\,(\eta^\bh,\xih,\theta^\bh):=\S'(\bus)[\bh]\,$
and $\,(\eta^\bk,\xik,\theta^\bk):=\S'(\bus)[\bk]$ are well-defined elements of the space $\,Z\times\tilde X\times Z$, where
$\,\S'(\bus)\,$ now denotes the extension of the Fr\'echet derivative introduced in Remark 2.5.  
We now claim that there is a constant $\widehat C>0$ that depends only on $\,R\,$ and the data, such that
\begin{equation}
\label{fiete10}
\|f_1\|_{L^2(Q)}+\|f_2\|_{L^2(Q)}\,\le\,\widehat C\,\|\bh\|_{L^2(Q)^2}\,\|\bk\|_{L^2(Q)^2}\,.
\end{equation}
Indeed, arguing as in the derivation of the estimates \eqref{fiete4}--\eqref{fiete7}, we obtain 
\begin{align*}
&\|f_1\|^2_{L^2(Q)}\,\le\,C\int_0^T \|\xik\|_4^2\,\left(\|\theta^\bh\|_4^2+\|\xih\|_4^2+\|\eta^\bh\|_4^2\right)\,ds\\
&\quad +\,C\int_0^T \|\xih\|_4^2\,\left(\|\theta^\bk\|_4^2+\|\xik\|_4^2+\|\eta^\bk\|_4^2\right)\,ds\\
&\quad +\,C\,\left(\|\s^*\|_{\Liq}^2+\|\vp^*\|_{\Liq}^2+\|\mu^*\|_{\Liq}^2+1\right)\int_0^T
\|\xik\|_4^2\,\|\xih\|_4^2\,ds\\[1mm]
&\,\le\,C\,\|\S'(\bus)[\bh]\|_{C^0([0,T];V)}^2\,\|\S'(\bus)[\bk]\|^2_{C^0([0,T];V)} \,\le\,C\,\|\bh\|_{L^2(Q)^2}^2
\,\|\bk\|^2_{L^2(Q)^2}\,,
\\[3mm]
&\|f_2\|_{L^2(Q)}^2\,\le\,C\int_0^T \|\xih\|_4^2\,\|\xik\|_4^2\,ds 
\,\le\,C\,\|\bh\|_{L^2(Q)^2}^2
\,\|\bk\|^2_{L^2(Q)^2}\,,
\end{align*}
which proves the claim. At this point, we can conclude from Lemma 2.4 that the system \Bilin\ has for every $\,\bh,\bk\in L^2(Q)^2\,$
a unique solution $(\nu,\psi,\rho)\in Z\times \widetilde X\times Z$. Moreover, we have, with a constant $K_6>0$ that depends
only on $\,R\,$ and the data, 
\begin{align}
\label{fiete11}
\|(\nu,\psi,\rho)\|_{Z\times \widetilde X\times Z}\,\le\,K_6\,\|\bh\|_{L^2(Q)^2}\,\|\bk\|_{L^2(Q)^2} \quad\forall\,
\bh,\bk\in L^2(Q)^2.
\end{align}

\vspace*{3mm}
\Brem
Similarly as in Remark 2.5, the operator $\,\S''(\bus)\in {\cal L}(\Uh,{\cal L}(\Uh,{\cal Y}))\,$ can be 
extended continuously to an element
of ${\cal L}(L^2(Q)^2,{\cal L}(L^2(Q)^2,Z\times\widetilde X\times Z))$ without changing its operator norm. Denoting 
the extended operator
still by $\,\S''(\bus)$, we see 
that the identity $\,\S''(\bus)[\bh,\bk]=(\nu,\psi,\rho)\,$ is also valid
for every $\,\bh,\bk \in L^2(Q)^2$, only that $\,(\nu,\psi,\rho)\in Z\times\widetilde X\times  Z$, in general. 
In addition, we have
\begin{equation}
\label{lsbound2}
\|\S''(\bu)[\bh,\bk]\|_{Z\times\widetilde X\times  Z}\,\le\,K_6\,\|\bh\|_{L^2(Q)^2}\,\|\bk\|_{L^2(Q)^2}\quad\mbox{for all 
\,$\bu\in\UR$\, and } \,\bh,\bk\in L^2(Q)^2\,.
\end{equation}
\Erem

\vspace*{2mm}
We conclude our preparatory work by showing a Lipschitz property for the extended operator $\,\S''\,$ that resembles \eqref{Lip1}.
\Blem
The mapping $\,\S'':{\cal U}\to {\cal L}(L^2(Q)^2,{\cal L}(L^2(Q)^2,Z\times\widetilde X\times Z))$, $\bu\mapsto \S''(\bu)$, is Lipschitz continuous in the
following sense: there is a constant $K_7>0$, which depends only on $\,R\,$ and the data, such that, for all controls $\,\bu_1,\bu_2
\in\UR\,$ and all increments $\bh,\bk\in L^2(Q)^2$,
\begin{align}
\label{Lip2}
\|\left(\S''(\bu_1)-\S''(\bu_2)\right)[\bh,\bk]\|_{\cal Z}\,\le\,K_7\,\|\bu_1-\bu_2\|_{L^2(Q)^2}\,\|\bh\|_{L^2(Q)^2}\,
\|\bk\|_{L^2(Q)^2}\,.
\end{align}
\Elem
\Bdim
We put $\,(\mu_i,\vp_i,\s_i):=\S(\bu_i)$, $\,(\eta_i^\bh,\xi_i^\bh,\theta_i^\bh):=\S'(\bu_i)[\bh]$, 
$\,(\eta_i^\bk,\xik,\theta_i^\bk):=\S'(\bu_i)[\bk]$, $(\nu_i,\psi_i,\rho_i):=\S''(\bu_i)[\bh,\bk]$, for $\,i=1,2$, as well as
\begin{align*}
&\bu:=\bu_1-\bu_2, \quad \mu:=\mu_1-\mu_2,\quad \vp:=\vp_1-\vp_2,\quad \s:=\s_1-\s_2,\\
&\eta^\bh:=\eta_1^\bh-\eta_2^\bh,\quad\xih:=\xih_1-\xih_2,\quad \theta^\bh:=\theta_1^\bh-\theta_2^\bh,\\
&\eta^\bk:=\eta_1^\bk-\eta_2^\bk,\quad\xik:=\xik_1-\xik_2,\quad \theta^\bk:=\theta_1^\bk-\theta_2^\bk,\\
&\nu:=\nu_1-\nu_2, \quad \psi:=\psi_1-\psi_2, \quad \rho:=\rho_1-\rho_2.
\end{align*} 
Then it follows from \eqref{stabu} and \eqref{Lip1} that
\begin{align}
\label{lippe0}
&\|(\mu,\vp,\s)\|_{{\cal Z}}\,\le\,C\,\|\bu\|_{L^2(Q)^2}\,, \quad \|(\eta^\bh,\xih,\theta^\bh)\|_{\cal Z}
\,\le\,C\,\|\bu\|_{L^2(Q)^2} \,\|\bh\|_{L^2(Q)^2}\,,        \nonumber\\[1mm]
&\|(\eta^\bk,\xik,\theta^\bk)\|_{\cal Z}
\,\le\,C\,\|\bu\|_{L^2(Q)^2} \,\|\bk\|_{L^2(Q)^2}\,.
\end{align}
We also recall the estimates \eqref{lsbound1} and \eqref{lsbound2}. Moreover, $\,(\nu,\psi,\rho)\,$ solves the problem
\begin{align}
\label{lippe1}
&\alpha \dt\nu+\dt\psi-\Delta\nu\,=\,P(\vp_1)(\rho-\chi\psi-\nu)+P'(\vp_1)(\s_1+\chi(1-\vp_1)-\mu_1)\psi\,+\,g_1\,,\\
\label{lippe2}
&\beta\dt\psi-\Delta\psi=\chi\rho-F''(\vp_1)\psi\,+\,g_2\,,\\
\label{lippe3}
&\dt\rho-\Delta\rho+\chi\Delta\psi=-P(\vp_1)(\rho-\chi\psi-\nu)-P'(\vp_1)(\s_1+\chi(1-\vp_1)-\mu_1)\psi\,+\,g_3\,,\\
\label{lippe4}
&\dn\nu=\dn\psi=\dn\rho=0\,,\\
\label{lippe5}
&\nu(0)=\psi(0)=\rho(0)=0\,,
\end{align}
which is again of the form \eqref{aux1}--\eqref{aux5} with $\lambda_1=\lambda_3=1$ and $\lambda_2=\lambda_4=0$, where
\begin{align}
\label{fritz1}
&g_1:=-g_3:=((P(\vp_1)-P(\vp_2))(\rho_2-\chi\psi_2-\nu_2)\,+\,P'(\vp_1)(\s-\chi\vp-\mu)\psi_2\nonumber\\
&+\,(P'(\vp_1)-P'(\vp_2))(\s_2+\chi(1-\vp_2)-\mu_2)\psi_2\,+\,(P'(\vp_1)-P'(\vp_2)) \,\xi_1^\bk\, (\theta_1^\bh-\chi\xi_1^\bh-\eta_1^\bh)
\nonumber\\
&+\,P'(\vp_2)\,\xik\,(\theta_1^\bh-\chi\xi_1^\bh-\eta_1^\bh)\,+\,P'(\vp_2)\,\xi_2^\bk\,(\theta^\bh-\chi\xi^\bh-\eta^\bh)
\nonumber\\
&+\,(P'(\vp_1)-P'(\vp_2)) \,\xi_1^\bh\, (\theta_1^\bk-\chi\xi_1^\bk-\eta_1^\bk)\,+\,
P'(\vp_2)\,\xih\,(\theta_1^\bk-\chi\xi_1^\bk-\eta_1^\bk)
\nonumber\\
&+\,P'(\vp_2)\,\xi_2^\bh\,(\theta^\bk-\chi\xi^\bk-\eta^\bk)\,+\,
(P''(\vp_1)-P''(\vp_2))\,\xi_1^\bh\,\xi_1^\bk\,(\s_1+\chi(1-\vp_1)-\mu_1)
\nonumber\\
&+\,P''(\vp_2)\,\xik\,\xi_1^\bh\,(\s_1+\chi(1-\vp_1)-\mu_1)\,+\,
P''(\vp_2)\,\xi_2^\bk\,\xih\,(\s_1+\chi(1-\vp_1)-\mu_1)
\nonumber\\
&+\,P''(\vp_2)\,\xi_2^\bk\,\xi_2^\bh\,(\s-\chi\vp-\mu)\,\,=:\,\,\sum_{i=1}^{13} B_i\,,
\\[3mm]
\label{fritz2}
&g_2:= -(F''(\vp_1)-F''(\vp_2))\psi_2-\left(F^{(3)}(\vp_1)-F^{(3)}(\vp_2)\right)\,\xi_1^\bh\,\xi_1^\bk\nonumber\\
&-F^{(3)}(\vp_2)\left(\xih\,\xi_1^\bk \,+\,\xi_2^\bh\,\xik\right)\,,
\end{align} 
where $\,B_i\,$ denotes the $i$th summand on the \rhs\ of \eqref{fritz1}.

At this point, we infer from the proof of \cite[Lem.~4.1]{CSS2} that the assertion follows once we can show that
$$
\sum_{i=1}^{13}\,\|B_i\|_{L^2(Q)}\,+\,\|g_2\|_{L^2(Q)}\,\le\,C\,\|\bu\|_{L^2(Q)^2}\,\|\bh\|_{L^2(Q)^2}\,
\|\bk\|_{L^2(Q)^2}\,.
$$
We only show the corresponding estimate for the terms $\,B_1,B_4,B_{11}\,$ and leave the others to the interested reader.
In the following, we make use of the mean value theorem, H\"older's inequality, the continuity of the embeddings 
$\,V\subset L^6(\Omega) \subset L^4(\Omega)$, and the global estimates \eqref{ssbound1}, \eqref{ssbound2},
\eqref{lsbound1}, and \eqref{lsbound2}.
We have
\begin{align*}
&\|B_1\|_{L^2(Q)}^2\,\le\,C\int_0^T\|\vp\|_4^2\left(\|\rho_2\|_4^2+\|\psi_2\|_4^2+\|\nu_2\|_4^2\right) ds\nonumber\\
&\quad\le\,C\,\|\vp\|^2_{C^0([0,T];V)}\,\|\S''(\bu_2)[\bh,\bk]\|^2_{C^0([0,T];V)^3}
\,\le\,C\,\|\bu\|_{L^2(Q)^2}^2\,\|\bh\|_{L^2(Q)^2}^2\,
\|\bk\|_{L^2(Q)^2}^2\,,\\[2mm]
&\|B_4\|^2_{L^2(Q)}\,\le\,C\int_0^T \|\vp\|_6^2\,\|\xi_1^\bk\|_6^2\left(\|\theta_1^\bh\|_6^2+\|\xi_1^\bh\|_6^2+\|\eta_1^\bh\|_6^2
\right) ds \\
&\quad \le\,C\,\|\vp\|^2_{C^0([0,T];V)}\,\|\xi_1^\bk\|_Z^2\,\|\S'(\bu_1)[\bh]\|^2_{\cal Z}
\,\le\,C\,\|\bu\|_{L^2(Q)^2}^2\,\|\bh\|_{L^2(Q)^2}^2\,\|\bk\|_{L^2(Q)^2}^2\,,
\end{align*}
as well as 
\begin{align*}
&\|B_{11}\|_{L^2(Q)}^2\,\le\,C\left(\|\s_1\|^2_{\Liq} + \|\vp_1\|^2_{\Liq} + \|\mu_1\|^2_{\Liq} + 1\right)
\int_0^T \|\vp\|_6^2\,\|\xi_1^\bh\|_6^2\,\|\xi_1^\bk\|_6^2\,ds
\\
&\le\,C\,\|\vp\|^2_{C^0([0,T];V)}\,\|\xi_1^\bh\|^2_{C^0([0,T];V)}\,\|\xi_1^\bk\|^2_{C^0([0,T];V)}
\,\le\,C\,\|\bu\|_{L^2(Q)^2}^2\,\|\bh\|_{L^2(Q)^2}^2\,\|\bk\|_{L^2(Q)^2}^2\,.
\end{align*}
The assertion of the lemma is thus proved.
\Edim


\section{The optimal control problem}

\setcounter{equation}{0}

We now begin to investigate the control problem {\bf (CP)}. In addition to {\bf (A1)}--{\bf (A6)}, we make the following assumptions:
\begin{description}
\item[(C1)] \,\,\,The constants $b_1,b_2$ are nonnegative, while $b_3,\kappa$ are positive.
\item[(C2)] \,\,\,It holds $\widehat\vp_\Omega\in \Huno$ and $\widehat\vp_Q\in L^2(Q)$.
\item[(C3)] \,\,\,$g:L^2(Q)^2\to\erre$ is nonnegative, continuous and convex on $L^2(Q)^2$.
\end{description}

\vspace{2mm}\noindent
Observe that {\bf (C3)} implies that $\,g\,$ is weakly sequentially lower 
semicontinuous on $L^2(Q)^2$. Moreover, denoting in the following by $\,\partial\,$ 
the subdifferential mapping in $L^2(Q)^2$, it follows from standard convex analysis that 
$\,\partial g\,$ is defined on the entire space $L^2(Q)^2$ and is a maximal monotone 
operator. In addition, 
the mapping
$((\mu,\varphi,\sigma),\bu)\mapsto J((\mu,\varphi,\sigma),\bu)$ 
defined by the cost functional \eqref{cost} is obviously
 continuous and convex 
(and thus weakly sequentially lower semicontinuous) on the
space $\bigl(L^2(Q)\times C^0([0,T];\Ldue)\times L^2(Q)\bigr)
\times L^2(Q)^2$. From a standard argument (which needs no repetition here) it then follows that
the problem {\bf (CP)} has a solution.

In the following, we often denote by $\bus=(u_1^*,u_2^*)\in \uad$ a local minimizer in the sense of $\Uh$
and by $\,(\mu^*,\vp^*,\s^*)=\S(\bus)\,$ the associated state. The corresponding  adjoint state variables 
solve the adjoint system, which is given by the backward-in-time parabolic system 
\begin{align}
\label{adj1}
	&- \dt p - \beta \dt q - \Delta q + \chi \Delta r  + F''(\vp^*)q -P'(\vp^*)(\s^*+\chi(1-\vp^*)-\mu^*)(p-r)\nonumber\\
	&\quad + \chi P(\vp^*)(p-r)
	=b_1 (\vp^* - \widehat \vp_Q) \quad\mbox{in }\,Q\,,\\[1mm]
	\label{adj2}
	&-\a \dt p-\Delta p -q + P(\vp^*)(p-r)=0 \quad\mbox{in }\,Q\,,\\[1mm]
	\label{adj3}
	& - \dt r -\Delta r  - \chi q - P(\vp^*)(p-r) =0 \quad\mbox{in }\,Q\,,\\[1mm]
	\label{adj4}
	&\dn p=\dn q=\dn r=0 \quad\mbox{on }\,\Sigma\,,\\[1mm] 
	\label{adj5}
	&(p+\b q)(T)= b_2(\vp^*(T)-\hat \ph_\Omega),
	\quad \a p(T)= 0,\quad r(T)=0,\quad \mbox{in }\,\Omega\,.
\end{align}
\Accorpa\Adjsys {adj1} {adj5}
According to \cite[Thm.~5.2]{CSS2}, the adjoint system has a unique weak 
solution $(\,p,q,r\,)$ satisfying
\begin{align}
\label{regs}
&p+\beta q\in H^1(0,T;V^*),\\
\label{regp}
&	p  \in \H1 H \cap C^0([0,T];V)\cap \L2 {W_0} \cap L^\infty(Q),\\
\label{regq}
&	q \in \L\infty H \cap \L2 V, 	\\
\label{regr}
&	r\in \H1 H \cap C^0([0,T];V) \cap \L2 {W_0}\cap L^\infty(Q),
\end{align}
as well as
\begin{align}
\label{wadj1}
	& -\langle \dt (p +\b q),v \rangle
	+ \iO \nabla q \cdot \nabla v
	- \chi \iO \nabla r  \cdot \nabla v
	+ \iO F''(\vp^*)\, q\,v
	\nonumber
	\\
	& -\iO P'(\vp^*)(\s^*+\chi(1-\vp^*)-\mu^*)
		(p-r) v \, +\, \chi \iO P(\vp^*)(p-r)v = b_1 \iO  (\vp^* - \hat \ph_Q) v,
	\\
	\label{wadj2}
	& -{\a} \iO\dt p\, v
	+ \iO \nabla p \cdot \nabla v 
	-\iO q\,v
	 +\iO  P(\vp^*)\,(p-r)\,v=0,
	\\
	\label{wadj3}
	& 
	- \iO \dt r \,v
	+\iO \nabla r\cdot \nabla v
	- \chi \iO q\, v
	 -\iO  P(\vp^*)\,(p-r)\, v=0 ,
\end{align}
for every $v\in V$ and almost every $t \in (0,T)$, and
\begin{equation}\label{wadj4}
		(p+\b q)(T)=b_2(\vp^*(T)-\hat \ph_\Omega),
		\quad p(T)= 0,\quad r(T)=0\quad \mbox{a.e. in }\,\Omega\,.
\end{equation}
Moreover, it follows from the proof of \cite[Thm.~5.2]{CSS2} that there exists a constant $\,K_8>0$, which depends only
on $\,R\,$ and the data (but not on the special choice of $\bus\in\uad$), such that
\begin{align}
\label{boundad}
&\|p\|_{H^1(0,T;H) \cap C^0([0,T];V)\cap L^2(0,T;W_0)\cap \Liq}\,+\,\|q\|_{H^1(0,T;V^*) \cap L^\infty(0,T;H)\cap L^2(0,T;V)} \nonumber\\
&\quad+\,\|r\|_{H^1(0,T;H) \cap C^0([0,T];V)\cap L^2(0,T;W_0)\cap \Liq} \nonumber\\
&\le\,K_8 \,\left(\|\vp^*-\hat\vp_Q\|_{L^2(Q)}\,+\,\|\vp^*(T)-\hat\vp_\Omega\|_{V}\right)\,.
\end{align}


\subsection{First-order necessary optimality conditions}

In this section, we aim at deriving  associated first-order necessary optimality conditions for local minima of the
optimal control problem {\bf (CP)}. We assume that {\bf (A1)}--{\bf (A6)} and {\bf (C1)}--{\bf (C3)} are fulfilled and
define  the reduced cost functionals associated with the functionals  $\,J\,$ and $\,J_1\,$  introduced in \eqref{cost} by
\begin{equation}\label{reduced}
\widehat J(\bu) =  J(\S(\bu),{\bf u}), \quad \widehat J_1(\bu)= J_1(\S(\bu),\bu)\,.
\end{equation}
Since $\S$ is twice continuously Fr\'echet differentiable  from $\,\Uh\,$ into $\,{\cal Y}\,$ and $\,{\cal Y}\,$ is continuously 
embedded in $\,C^0([0,T];L^2(Q)^3)$,
$\S\,$ is also twice continuously Fr\'echet differentiable from $\,\Uh\,$ into $\,C^0([0,T];L^2(Q)^3)$. It thus follows 
from the chain rule that the smooth part $\,\widehat J_1\,$ of $\,\widehat J\,$ is  a twice continuously Fr\'echet 
differentiable mapping from $\Uh$ into $\erre$, where, for every $\bus=(u_1^*,u_2^*)\in\Uh$ and every ${\bf h}=(h_1,h_2)\in\Uh$, it holds 
with $(\mu^*,\vp^*,\s^*)=\S(\bus)$ that
\begin{align}
\label{DJ}
\widehat J'_1(\bus)[{\bf h}]\,&=\,b_1\iint_Q \xih(\vp^*-\hat\vp_Q)\,+\,b_2\iO \xih(T)(\vp^*(T)-\hat\vp_{\Omega})  \nonumber\\
&\quad + b_3 \iint_Q (u_1^* h_1 \,+\,u^*_2 h_2),
\end{align}
where $\,(\eta^\bh,\xih,\theta^\bh)=\S'(\bus)[\bf h]\,$ is the unique solution to the linearized system \eqref{aux1}--\eqref{aux5}, with
$\,\lambda_1=\lambda_2=1\,$ and $\,\lambda_3= \lambda_4=0$,
associated with ${\bf h}$.
\Brem
Observe that the \rhs\ of \eqref{DJ} is meaningful also for arguments ${\bf h}=(h_1,h_2)\in L^2(Q)^2$, where in 
this case $(\eta^\bh,\xih,\theta^\bh)=\S'(\bus)[{\bf h}]$ with the extension of the operator $\,\S'(\bus)$\, to $\,L^2(Q)^2$ introduced in Remark 2.5. Hence, by means of the
identity \eqref{DJ} we can extend the operator $\,\widehat J_1'(\bus)\in\Uh^*\,$ to $\,L^2(Q)^2$. 
The extended operator, which we again denote by
$\widehat J'_1(\bus)$, then becomes an element of $(L^2(Q)^2)^*$. In this way, expressions of the form $\widehat 
J_1'(\bus)[{\bf h}]$  have a proper meaning also for ${\bf h}\in L^2(Q)^2$.  
\Erem

In the following, we assume that $\bus = (u_1^*,u_2^*)\in\uad$ is a given locally optimal control for ${\bf (CP)}$ in the sense
of \,$\Uh$, that is, there is some $\varepsilon>0$ such that
\begin{equation}
\label{lomin}
\widehat J(\bu)\,\ge\,\widehat J(\bus)\quad\mbox{for all $\bu\in\uad$ satisfying }\,\|\bu-\bus\|_{\Uh}\,\le\,\varepsilon.
\end{equation}
Notice that any locally optimal control in the sense of $\,L^p(Q)^2\,$ with $1 \le p < \infty$
is also locally optimal in the sense of $\Uh$. 
Therefore, a result proved for locally optimal controls in the sense of  $\,\Uh\,$ is also valid for locally optimal controls 
in the sense of  $\,L^p(Q)^2$. It is of course also valid for  (globally) optimal controls.

Now, in the same way as in \cite{SpTr2}, we infer that then  the variational inequality
\begin{equation}
\label{var1}
{\widehat J}_1'(\bus)[\bu -  \bus] + \kappa\left(g(\bu) - g(\bus)\right) \ge 0 \quad \forall\, \bu \in \uad
\end{equation}
is satisfied. Moreover, 
denoting by the symbol $\,\partial\,$ the subdifferential mapping in $L^2(Q)^2$ (recall that $g$ is a convex 
continuous functional on  $L^2(Q)^2$), we conclude from \cite[Thm.~4.5]{SpTr2} that there is 
some $\bl^*=(\lambda^*_1,\lambda_2^*) \in \partial g(\bus)\subset L^2(Q)^2$ such that 
\begin{equation} \label{varineq1}
{\widehat J}_1'(\bus)[\bu - \bus] +  \iint_Q \kappa\,(\lambda_1^* (u_1-u_1^*)\,+\, \lambda_2^*(u_2-u^*_2)) \,\ge 0 \quad \forall \,
\bu=(u_1,u_2) \in \uad.
\end{equation}
As usual, we simplify the expression $\,{\widehat J}_1'(\bus)[\bu-\bus]\,$ in \eqref{varineq1}
by means of the adjoint state variables defined in \Adjsys. A standard calculation (see the proof of \cite[Thm.~5.4]{CSS2}) then leads 
to the following result.
\Bthm
\,\,{\rm (Necessary optimality condition)}  Suppose that {\bf (A1)}--{\bf (A6)} and {\bf (C1)}--{\bf (C3)}
are fulfilled. Moreover, 
 let $\bus=(u_1^*,u_2^*) \in \uad$ be a locally optimal control of  {\bf (CP)} in the sense of $\,\Uh\,$ 
with associated state $\,(\mu^*,\vp^*,\s^*)={\cal S}(\bus)$
and adjoint state $\,(p^*,q^*,r^*)$.
Then there exists some ${\boldsymbol \lambda}^*=(\lambda_1^*,\lambda_2^*)   \in \partial g(\bus)$ such that,
for all $\, \bu =(u_1,u_2)\in \uad$, 
\begin{align}
\label{varineq2}
&\iint_Q  (-\h\,p^*+\kappa\lambda^*_1+b_3 u^*_1 )(u_1-u^*_1) \,+\,\iint_Q
(r^*+\kappa\lambda_2^*+b_3u^*_2) (u_2-u_2^*) \, \ge \,0 \,.
\end{align}
\Ethm
\Brem We underline again that \eqref{varineq2} is also necessary for all globally optimal controls and all controls 
which are even locally optimal in the sense of $L^p(Q)\times L^p(\Sigma)$ with $p \ge 1$.
Observe also that the variational inequality \eqref{varineq2} is equivalent to two independent variational inequalities for $u_1^*$ and $u_2^*$ 
that have to hold simultaneously, namely,
\begin{eqnarray}
\label{varin1}
\iint_Q  \left(-\h\, p^*+\kappa\lambda^*_1+b_3 u_1^*\right)\left(u_1 - u_1^*\right) &\!\!\ge\!\!& 0 \quad \forall\, u_1 \in U_{\rm ad}^1, \\
\label{varin2}
\iint_Q  \left(r^*+\kappa\lambda^*_2 +b_3u_2^*\right)\left(u_2 - u_2^*\right) 
&\!\!\ge\!\!& 0 \quad \forall\, u_2 \in U_{\rm ad}^2, 
\end{eqnarray}
where 
\begin{align}
\label{uadi}
U_{\rm ad}^i\,& := \,\{ u_i \in L^\infty(Q): \underline u_i \, \le u_i \le \,\overline u_i \,\, \mbox{ a.e. in }\,  Q\},\quad i=1,2.
\end{align}
\Erem

\subsection{Sparsity of controls}

The convex function $\,g\,$ in the objective functional accounts for the sparsity of optimal controls, i.e., 
any locally optimal control can vanish in some region of the space-time cylinder $Q$. The form of this region depends on the particular choice of the functional $\,g\,$ which can differ in different situations.
The sparsity properties can be deduced from the variational inequalities \eqref{varin1} and \eqref{varin2} and the 
form of the subdifferential  $\partial g$. In this paper, we restrict our analysis to the case of {\em full sparsity} which is 
characterized by the functional (recall \eqref{cost}) 
\begin{equation}
\label{defg}
g(\bu)=g(u_1,u_2):=\iint_Q\left(|u_1|+|u_2|\right)\,.
\end{equation}
Other important choices leading to the {\em directional sparsity with respect to time} and the {\em directional sparsity with respect to space}
are not considered here. It is well known (see, e.g., \cite{ST}) that the subdifferential of $\,g\,$ is given by  
\begin{align} 
\label{dg}
&\partial g(\bu) = g(u_1,u_2)\nonumber\\[2mm]
&:=\left\{(\lambda_1,\lambda_2) \in L^2(Q)^2:\,
\lambda_i \in \left\{
\begin{array}{ll}
\{1\} & \mbox{ if } u_i > 0\\
{[-1,1]}& \mbox{ if } u_i = 0\\
\{-1\} & \mbox{ if } u_i < 0\\
\end{array}
\right\}
\mbox{ a.e. in }\,Q, \,\,\,i=1,2
\right\}\,.
\end{align}

The following sparsity result can be proved in exactly the same way as \cite[Thm.~4.9]{SpTr2}.

\Bthm {\rm (Full sparsity)} \,\,\,Suppose that the assumptions {\bf (A1)}--{\bf (A6)} and {\bf (C1)}--{\bf (C3)} are fulfilled, and assume that 
$\,\underline u_i<0 <\overline u_i$,  $i=1,2$. Let $\bus=(u_1^*,u_2^*)\in \uad$ be a locally optimal control in the sense 
of \,$\Uh$\, for the problem {\bf (CP)}
with the sparsity functional $\,g\,$ defined in \eqref{defg}, and with associated state $(\mu^*,\vp^*,\s^*)=\S(\bus)$ solving \Statesys\ and 
adjoint state $(p^*,q^*,r^*)$ solving \Adjsys. Then there exists some 
$(\lambda^*_1,\lambda_2^*)\in\partial g(\bus)$  such that \eqref{varin1}--\eqref{varin2} are satisfied.
In addition, we have that
\begin{eqnarray}
u_1^*(x,t) = 0 \quad &\Longleftrightarrow& \quad |-\h(x,t)p^*(x,t)| \le \kappa, \quad\mbox{for a.e. }\,(x,t)\in Q, 
\label{u1sparsity}
\\
\label{u2sparsity}
u_2^*(x,t) = 0 \quad &\Longleftrightarrow& \quad |r^*(x,t)| \le \kappa, \quad\mbox{for a.e. }\,(x,t)\in Q. 
\end{eqnarray}
Moreover, if\, $(p^*,q^*,r^*)$ and $(\lambda_1^*, \lambda_2^*)$ are given, then
$(u_1^*,u_2^*)$ is obtained from the projection formulas
\begin{align}
\label{pro1} 
u_1^*(x,t)& = \max\left\{\underline u_1, \min\left\{ \overline u_1, 
-{b_3}^{-1} \left(-\h\,p^*+ \kappa \,\lambda_1^*\right)(x,t)\right\}\right\}
\,\mbox{ for a.e. $(x,t)\in Q$} ,\\[1mm]
\label{pro2}
u_2^*(x,t) &=
 \max\left\{\underline u_2,\min\left\{ \overline u_2,-{b_3}^{-1} \left(r^*+\kappa\lambda_2^*\right)(x,t)\right\}\right\}
\,\mbox{ for a.e. $(x,t) \in Q$.}
\end{align}
\Ethm

The projection formulas above are standard conclusions from the variational inequalities \eqref{varin1}--\eqref{varin2}.
 
\subsection{Second-order sufficient optimality conditions} 
In this section, we establish the main results of this paper, using auxiliary results collected in the Appendix. 
We provide conditions that ensure local optimality of pairs \,$\bus=(u^*_1,u^*_2)$ obeying the first-order necessary optimality 
conditions of Theorem 3.2. Second-order sufficient optimality conditions are based on a condition of coercivity that 
is required to hold for the smooth part  $\,\widehat J_1\,$ 
of $\,\widehat {\cal J}\,$ in a certain critical cone. The nonsmooth part $\,g\,$ contributes to sufficiency by its convexity. In the following,
we generally assume that {\bf (A1)}--{\bf (A6)}, {\bf (C1)}--{\bf (C3)}, and the conditions
$\,\underline u_1<0<\overline u_1\,$ and $\,\underline u_2<0<\overline u_2\,$ are fulfilled.
Our analysis will follow closely the lines of \cite{SpTr2}, which in turn follows \cite{casas_ryll_troeltzsch2015}, where a second-order analysis was 
performed for sparse control of the FitzHugh--Nagumo system. In particular, we adapt the proof of 
\cite[Thm.~3.4]{casas_ryll_troeltzsch2015} to our setting of less regularity.

To this end, we fix  a pair of controls $\,\bus = (u_1^*,u_2^*)$ that satisfies the first-order necessary optimality conditions,
and we set $\,(\mu^*,\vp^*,\s^*)=\S(\bus)$. Then the cone
\[
C(\bus) = \{ (v_1,v_2) \in L^2(Q)^2 \,\text{ satisfying the sign conditions \eqref{sign} a.e. in $Q$}\},
\]
where
\begin{equation} \label{sign}
v_i(x,t) \left\{
\begin{array}{l}
\ge 0 \,\, \text{ if }\,\, u^*_i(x,t) = \underline u_i\\[1mm] 
\le 0 \,\, \text{ if }\,\, u^*_i(x,t) = \overline u_i 
\end{array}
\right. \,,\quad i=1,2\,,
\end{equation}
is called the {\em cone of feasible directions}, which is a convex and closed subset of $L^2(Q)^2$.
 We also need the directional derivative of $\,g\,$ at $\bu\in L^2(Q)^2$ in the direction 
$\bv=(v_1,v_2)\in L^2(Q)^2$, which is given by
\begin{equation}
\label{g'}
g'(\bu,\bv) = \lim_{\tau \searrow 0} \frac{1}{ \tau}(g(\bu+\tau \bv)-g(\bu))\,.
\end{equation}
Following the definition of the critical cone in \cite[Sect. 3.1]{casas_ryll_troeltzsch2015}, we define
\begin{equation}
\label{critcone}
C_{\bus} = \{\bv \in C(\bus): D\widehat{J}(\bus)[\bv] + \kappa g'(\bus,\bv) = 0\}\,,
\end{equation}
which is also a closed and convex subset of $L^2(Q)^2$. According to \cite[Sect. 3.1]{casas_ryll_troeltzsch2015}, 
it consists of all $\bv=(v_1,v_2)\in C(\bus)$ satisfying  
\begin{align} \label{pointwise1}
&v_1(x,t) \left\{
\begin{array}{l}
= 0 \,\,\mbox{ if } \,\,|-\h(x,t)p^*(x,t) + b_3 u_1^*(x,t)| \not= \kappa\\[1mm]
\ge 0 \,\,\mbox{ if }\,\, u_1^*(x,t) = \underline u_1\,\, \mbox{ or } \,\,(-\h(x,t)p^*(x,t) = -\kappa \,\,\mbox{ and }\,\, u_1^*(x,t) = 0)\\[1mm]
\le 0 \,\,\mbox{ if }\,\, u_1^*(x,t) = \overline u_1\,\, \mbox{ or }\,\, (-\h(x,t)p^*(x,t) = \kappa \,\,\mbox{ and } \,\,u_1^*(x,t) = 0)
\end{array}
\right.,\\[1mm]
\label{pointwise2}
&v_2(x,t) \left\{
\begin{array}{l}
= 0 \,\,\mbox{ if } \,\,|r^*(x,t) + b_3 u_2^*(x,t)| \not= \kappa\\[1mm]
\ge 0 \,\,\mbox{ if }\,\, u_2^*(x,t) = \underline u_2\,\, \mbox{ or } \,\,(r^*(x,t) = -\kappa \,\,\mbox{ and }\,\, u_2^*(x,t) = 0)\\[1mm]
\le 0 \,\,\mbox{ if }\,\, u_2^*(x,t) = \overline u_2\,\, \mbox{ or }\,\, (r^*(x,t) = \kappa \,\,\mbox{ and } \,\,u_2^*(x,t) = 0)
\end{array}
\right. .
\end{align}

\Brem Let us compare the first condition in \eqref{pointwise1} with the situation in the differentiable control problem without sparsity terms obtained for $\kappa=0$.
Then this condition leads to the requirement that \,\,$v_1(x,t) = 0 \mbox{ if } \ |-\h(x,t)p^*(x,t) + b_3 u^*_1(x,t)|> 0$,
or, since $\kappa=0$,
\begin{equation}  \label{simplecone}
v_1(x,t) = 0 \mbox{ if } \ |-\h(x,t)p^*(x,t) + \kappa \lambda^*_1(x,t) + b_3 u^*_1(x,t)|> 0. 
\end{equation}
An analogous condition results for $\,v_2$.

One might be tempted to define the critical cone using \eqref{simplecone} and its counterpart for $\,v_2\,$ 
also in the case $\kappa > 0$.
This, however, is not a good idea, because it leads to a critical cone that is larger than needed, in general. 
As an example, we mention the particular case when the control
$\bus = \mathbf{0}$ satisfies the first-order necessary optimality conditions and when $\,\,|-\h\,p^*| < \kappa\,\,$
and \,\,$|r^*|<\kappa$\,\, hold a.e. in $Q$. 
Then the upper relation of \eqref{pointwise1}, and its counterpart for $v_2$, lead to $C_{\bus}=\{\mathbf{0}\}$,
the smallest possible critical cone.

However, thanks to $u^*_1=0$, the variational inequality \eqref{varin1} implies that $-\h\,p^* + \kappa \lambda^*_1 + b_3 u^*_1=0$ 
a.e. in $Q$, i.e., the condition \,$|-\h(x,t)p^*(x,t) + \kappa \lambda^*_1(x,t) + b_3 u^*_1(x,t)|> 0\,$ can only be satisfied on a set of measure zero. Moreover,
also the sign conditions \eqref{sign} do not restrict the critical cone. Hence, the largest possible critical cone $C_{\bus} = L^2(Q)^2$ would be obtained, 
provided that analogous conditions hold for $u_2^*$ and $r^*$ in $\,Q$.

In this example, the quadratic growth condition \eqref{growth} below is valid for the choice \eqref{critcone} as critical cone even 
without assuming the coercivity condition 
\eqref{coerc} below (here the so-called first-order sufficient conditions apply), while the use of a cone based on 
\eqref{simplecone} leads to postulating \eqref{coerc} on the whole space $L^2(Q)^2$ for the quadratic growth 
condition to be valid. This shows that the choice of \eqref{critcone} as critical cone is essentially better than 
of one based on 
\eqref{simplecone}. 
\Erem

At this point, we derive an explicit expression for $\,\widehat J_1''(\bu)[\bv,\mathbf{w}]\,$ for arbitrary $\,\bu=(u_1,u_2),
\bv=(v_1,v_2), \mathbf{w}=(w_1,w_2)\in\Uh$. In the following, we argue similarly as in 
\cite[Sect.~5.7]{Fredibuch} (see also \cite[Sect.~6]{CSS2}). At first, we readily infer that, 
for every $((\m,\ph,\s), \bu)\in (\C0 H)^3\times \Uh$ 
and ${\bf v}=(v_1,v_2,v_3),{\bf w}=(w_1,w_2,w_3)$ such that
$ ({\bf v}, \bh),({\bf w}, \bk)\in (\C0 H)^3\times \Uh$, we have 
\begin{align}
	\label{D2:0}
	J_1''((\m,\ph,\s), \bu)[({\bf v},\bh),({\bf w},\bk)]
	= b_1 \iint_Q  v_2 w_2 \,+\,b_2 \iO v_2(T) w_2(T)
	\,+\,b_3 \iint_Q \bh\cdot \bk,
\end{align}
where the dot denotes the euclidean scalar product in $\erre^2$. For the second-order derivative of the reduced cost functional $\widehat J_1$
at a fixed control $\bus$ we then find with $\,(\mu^*,\vp^*,\s^*)=\S(\bus)$\, that
\begin{align}
	\non \widehat J_1''(\bus)[\bh,\bk] & = D_{(\m,\ph,\s)}J_1((\mu^*,\vp^*,\s^*), \bus)[(\nu,\psi,\rho)]\\
	\label{D2:1}
	&\quad + \,J_1'' ((\bm,\bph,\bs), \bus)[((\eta^\bh,\xih,\th^\bh),\bh),((\eta^\bk,\xik,\th^\bk),\bk)],
\end{align}
where $(\eta^\bh,\xih,\th^\bh)$, $(\eta^\bk,\xik,\th^\bk)$, and $(\nu,\psi,\rho)$ stand for the unique corresponding solutions to the linearized 
system associated with $\bh$ and $\bk$}, and to the bilinearized system, respectively.
From the definition of the cost functional \eqref{cost} we readily infer that
\begin{align}
	\label{D2:2} 
	D_{(\m,\ph,\s)} J_1((\bm,\bph,\bs),\bus)[(\nu,\psi,\rho)] 
	= b_1   \iint_Q (\bph - \hat \ph_Q) \psi
	\,+\, b_2\iO (\bph(T) - \hat \ph_{\Omega})\psi(T).
\end{align}
We now claim that, with the associated adjoint state $\,(p^*,q^*,r^*)$,
\begin{align}
&b_1 \iint_Q (\bph-\hat \ph_Q)\psi 
 	\,+ b_2\iO (\bph(T) - \hat \ph_{\Omega})\psi(T)\nonumber\\
& = \iint_Q \big[ P'(\bph)\left(\xik(\th^\bh-\chi\xih-\eta^\bh)+\xih(\th^\bk-\chi\xik-\eta^\bk)\right)(p^*-r^*)\nonumber\\
&\hspace*{10mm} + P''(\bph)\xik\xih(\bs +\chi (1-\bph)-\bm)(p^*-r^*) \,-\, F^{(3)}(\bph)\xih\xik q^* \big].
\end{align}
To prove this claim, we multiply \eqref{bilin1} by $\,p^*$, \eqref{bilin2} by $\,q^*$, \eqref{bilin3} by $\,r^*$, 
add the resulting equalities, and integrate over $\,Q$, to obtain that
\begin{align*}
	0& = \iint_Q p^* \,\Big[\alpha\dt\nu+\dt\psi-\Delta\nu - P(\vp^*)(\rho-\chi\psi-\nu)\nonumber\\
	&\hspace*{19mm} -P'(\vp^*)(\s^*+\chi(1-\vp^*)-\mu^*)\psi \,- \,f_1 \Big]
	\\ 
	& \quad
	+ \iint_Q q^* \,\Big[\beta\dt\psi-\Delta\psi-\nu-\chi \rho + F''(\bvp)\psi + F^{(3)}(\bph) \xih \xik \Big]
	\\ 
	& \quad
	+ \iint_Q r^*\, \Big[\dt\rho-\Delta\rho+\chi\Delta\psi + P(\vp^*)(\rho-\chi\psi-\nu)\nonumber\\
	&\hspace*{21mm} +P'(\vp^*)(\s^*+\chi(1-\vp^*)-\mu^*)\psi \,	+\,f_1 \Big]
	\end{align*} 
with the function $\,f_1\,$ defined in \eqref{fiete8}. Then, we integrate by parts and make use of the initial and terminal conditions 
\eqref{bilin5} and \eqref{wadj4} to find that
\begin{align*}
0& = \iint_Q \nu \,\Big[-\a \dt p^*-\Delta p^* -q^* + P(\bph)(p^*-r^*)\Big]
\\ 
&\quad
+ \int_0^T \langle -\dt(p^*+\b q^*)(t),\psi(t)\rangle \,dt \,+\,b_2\iO(\vp^*(T)-\widehat\vp_\Omega)\psi(T)\\
&\quad +
\iint_Q \psi\,\Big [- \Delta q^* + \chi \Delta r^*  + F''(\bph) q^* + \chi P(\bph)(p^*-r^*)\\
& \hspace*{23mm}
	-P'(\bph)(\bs+\chi(1-\bph)-\bm)(p^*-r^*)  \Big]
\\ 
& \quad
+ \iint_Q \rho\,\Big [- \dt r^* -\Delta r^*  - \chi q^* - P(\bph)(p^*-r^*)\Big]
\\ 
& \quad
+ \iint_Q \Big[- P'(\bph)\left(\xik(\th^\bh-\chi\xih-\eta^\bh)+ \xih(\th^\bk-\chi\xik-\eta^\bk)\right)(p^*-r^*)
\\ 	
&\qquad\qquad\quad - P''(\bph)\xik\xih(\bs +\chi (1-\bph)-\bm)(p^*-r^*)
 	+ F^{(3)}(\bph)\xih\xik q^* \Big],
\end{align*}
whence the claim follows, since $\,(p^*,q^*,r^*)\,$ solves the adjoint system \eqref{wadj1}--\eqref{wadj4}.
From this characterization, along with \eqref{D2:1} and \eqref{D2:2}, we conclude that
\begin{align}	
& \widehat J_1''(\bus)[\bh,\bk] \,=\,b_1\iint_Q \xih\,\xik \,+\,b_2\iO \xih(T)  \xik(T)
\,+\,b_3\iint_Q \bh\cdot\bk\nonumber\\
& +\iint_Q \Big[ P'(\bph)\left(\xik(\th^\bh-\chi\xih -\eta^\bh)+\xih(\th^\bk-\chi\xik-\eta^\bk)\right)(p^*-r^*)\nonumber\\
&\hspace*{15mm}	+ P''(\bph)(\bs +\chi (1-\bph)-\bm)(p^*-r^*)\xih\xik \,-\,F^{(3)}(\bph) \xih\xik q^*\Big]	\,.
\label{walter1}
\end{align}

Observe that the expression on the right-hand side of \eqref{walter1} is meaningful also for increments $\,\bh,\bk\in L^2(Q)^2$. 
Indeed, in this case the expressions $\,(\eta^\bh,\xih,\theta^\bh)=\S'(\bus)[\bh]$, $\,(\eta^\bk,\xik,\theta^\bk)=\S'(\bus)[\bk]$,
and $\,(\nu,\psi,\rho)=\S''(\bus)[\bh,\bk]\,$ have an interpretation in the sense of the extended operators $\,\S'(\bus)\,$ and
$\,\S''(\bus)\,$ introduced in Remark 2.5 and Remark 2.7. Therefore, the operator $\,\widehat J_1''(\bus)\,$ can be extended by
the identity \eqref{walter1} to the space $L^2(Q)^2\times L^2(Q)^2$. This extension,
which will still be denoted by $\,\widehat J_1''(\bus)$, will be frequently used in the following.
We now show that it is continuous. Indeed, we claim 
that for all $\bh,\bk\in L^2(Q)^2$ it holds 
\begin{align}
\label{walter2}
\left|\widehat J_1''(\bus)[\bh,\bk]\right|\,\le\,\widehat C\,\|\bh\|_{L^2(Q)^2}\,\|\bk\|_{L^2(Q)^2}\,,
\end{align}
where the constant $\widehat C>0$ is independent of the choice of $\,\bus \in {\cal U}_R$. Obviously, only the last integral on the \rhs\
of \eqref{walter1} needs some treatment, and we estimate just its third summand, leaving the others as an exercise to the reader. We have,
by virtue of H\"older's inequality, the continuity of the embedding $V\subset L^4(\Omega)$, and the global bounds \eqref{ssbound3}, \eqref{lsbound1}, and \eqref{boundad}, 
\begin{align*}
&\Big| \iint_Q F^{(3)}(\bph) \xih\xik q^* \Big|\,\le\,C\int_0^T \|\xih\|_4\,\|\xik\|_4\,\|q^*\|_2\,dt\\[2mm]
&\le \,C\,\|\xih\|_{C^0([0,T];V)}\,\|\xik\|_{C^0([0,T];V)}\,\|q^*\|_{L^\infty(0,T;H)}\,\le\,
C\,\|\bh\|_{L^2(Q)^2}\,\|\bk\|_{L^2(Q)^2}\,, 
\end{align*}
as asserted. 

In the following, we will employ the following coercivity condition:

\begin{equation} \label{coerc}
\widehat J_1''(\bus)[\bv,\bv] >0 \quad \forall \, \bv \in C_{\bus} \setminus \{\mathbf{0}\}\,.
\end{equation}

\noindent Condition \eqref{coerc} is a direct extension of associated conditions that are standard in finite-dimensional 
nonlinear optimization. In the optimal control of partial differential equation, it was first used in \cite{casas_troeltzsch2012}.
As in  \cite[Thm 3.3]{casas_ryll_troeltzsch2015} or \cite{casas_troeltzsch2012}, it can be shown that \eqref{coerc} is equivalent to the existence of a 
constant $\,\delta > 0\,$ such that $\,\widehat J_1''(\bus)[\bv,\bv] \ge \delta \, \|\bv\|^2_{L^2(Q)^2}\,$ for all $\,\bv \in C_{\bus}$.

We have the following result.
\Bthm \,\,{\rm (Second-order sufficient condition)} \,\,Suppose that {\bf (A1)}--{\bf (A6)} and {\bf (C1)}--{\bf (C3)}
are fulfilled and that
$\,\underline u_i<0<\overline u_i$, $\,i=1,2$. Moreover, let $\bus=(u^*_1,u_2^*) \in \uad$, together with the associated state 
$(\mu^*,\vp^*,\s^*)=\S(\bus)$ and the adjoint state $(p^*,q^*,r^*)$, fulfill the first-order necessary optimality conditions 
of Theorem 3.2. If, in addition,  $\bus$ satisfies the coercivity condition \eqref{coerc}, 
then there exist constants $\,\varepsilon > 0\,$ and $\,\tau > 0\,$ such that the quadratic growth condition
\begin{equation} \label{growth}
\widehat J(\bu) \ge \widehat J(\bus) + \tau \, \|\bu-\bus\|^2_{L^2(Q)^2} 
\end{equation}
holds for all $\,\bu \in \uad\,$ with $ \|\bu-\bus\|_{L^2(Q)^2}  < \varepsilon$. Consequently,
$\,\bus\,$ is a locally optimal control in the sense of $\,L^2(Q)^2$.
\Ethm
\begin{proof} The proof follows that of \cite[Thm.~3.4]{casas_ryll_troeltzsch2015}. 
We argue by contradiction, assuming that the claim of the theorem is not true. Then there exists a sequence of controls 
$\{\bu_k\}\subset \uad$ such that, for all $k\in\enne$, 
\begin{equation}
\label{contrary}
\|\bu_k-\bus\|_{L^2(Q)^2} < \frac{1}{k} \quad \mbox{ while } \quad \widehat J(\bu_k)<  \widehat J(\bus) 
+ \frac{1}{2k} \|\bu_k-\bus\|_{L^2(Q)^2}^2\,.
\end{equation}
Noting that $\bu_k\not=\bus$ for all $k\in\enne$, we define 
$$
r_k = \|\bu_k-\bus\|_{L^2(Q)^2}
  \quad \mbox{ and } \quad \bv_k = \frac{1}{r_k}(\bu_k-\bus)\,.
$$                                               
Then $\|\bv_k\|_{L^2(Q)^2}=1$ and, possibly after selecting a subsequence, we can assume that 
\[
\bv_k \to \bv \, \mbox{ weakly in }\,  L^2(Q)^2
\]
for some $\bv\in L^2(Q)^2$. As in \cite{casas_ryll_troeltzsch2015}, the proof is split into three parts. 

(i) $\bv \in C_{\bus}$: Obviously, each $\bv_k$ obeys the sign conditions \eqref{sign} and thus belongs to $C(\bus)$. 
Since $C(\bus)$ is convex and closed in $L^2(Q)^2$, it follows that $\bv\in C(\bus)$. We now claim that 
\Beq
\label{Paulchen}
\widehat J_1'(\bus)[\bv] + \kappa g'(\bus,\bv) = 0.
\end{equation}
Notice that by Remark 3.1 the expression 
$\,\widehat{J}_1'(\bus)[\bv]\,$ is well defined. For every $r \in (0,1)$ and all $\bv=(v_1,v_2),\, \bu=(u_1,u_2) \in L^2(Q)^2$, 
we infer from the convexity of $\,g\,$ that
\begin{align}
g(\bv)-g(\bu) &\,\ge\, \frac{g(\bu + r (\bv-\bu))-g(\bu)}{r} \,\ge\, g'(\bu,\bv-\bu)\nonumber \\
&\,=\, \max_{(\lambda_1,\lambda_2)\in \partial g(\bu)}\,
\iint_Q \Big(\lambda_1 (v_1-u_1) + \lambda_2 (v_2-u_2)\Big).
\label{directionalder}
\end{align}
In particular, with $\bu_k=(u_{k_1},u_{k_2})$,
\begin{align}
&\widehat{J}_1'(\bus)[\bv] +  \kappa g'(\bus,\bv)\,\ge\, \widehat{J}_1'(\bus)[\bv] + \iint_Q \kappa\big(\lambda^*_1 v_1 \,+\, 
\lambda_2^* v_2\big) \nonumber\\
&= \iint_Q  \big((-\h p^* + b_3 u^*_1 +\kappa\lambda_1^*) v_1 \,+ \, (r^* + b_3 u_2^*+\kappa\lambda_2^*) v_2\big) \nonumber\\
&= \lim_{k \to \infty} \frac{1}{r_k}  \iint_Q \big( (-\h p^* + b_3 u^*_1 + \kappa \lambda^*_1) (u_{k_1}-u^*_1)
\,+\,(r^* + b_3 u_2^* + \kappa \lambda_2^*) (u_{k_2}-u_2^*)\big)\nonumber\\
&\ge 0\,,
\end{align}
by the variational inequality \eqref{varineq2}. Next, we prove the converse inequality. By \eqref{contrary}, we have
\[
\widehat{J}_1(\bu_k)-\widehat{J}_1(\bus) + \kappa\left(g(\bu_k)-g(\bus)\right) < \frac{1}{2k} r_k^2\,,
\]
whence, owing to the mean value theorem, and since $\bu_k = \bus+ r_k \bv_k$, we get
\[
\widehat{J}_1(\bus) + r_k \widehat{J}_1'(\bus + \theta_k r_k \bv_k)[\bv_k] + \kappa g(\bus + r_k \bv_k)
< \widehat{J}_1(\bus) + \kappa g(\bus) +  \frac{1}{2k} r_k^2
\]
with some $0 < \theta_k < 1$. From \eqref{directionalder}, we obtain $\kappa(g(\bus + r_k \bv_k)-g(\bus))\ge  
\kappa g'(\bus,r_k \bv_k)$, and thus
\[
r_k \widehat{J}_1'(\bus + \theta_k r_k \bv_k)[\bv_k] + r_k \kappa g'(\bus,\bv_k) < \frac{ r_k^2}{2k}\,.
\]
We divide this inequality by $r_k$ and pass to the limit $k \to \infty$. Here,  we invoke Corollary \ref{coroll3} of the Appendix, and
we use that $\,g'(\bus,\bv_k)\to g'(\bus,\bv)$. We then obtain the desired converse inequality
\[
\widehat{J}_1'(\bus)[\bv] + \kappa g'(\bus,\bv) \le 0\,,
\]
which completes the proof of (i).

(ii) $\bv = {\bf 0}$: We again invoke \eqref{contrary}, now  performing a second-order Taylor expansion on the left-hand side,
\begin{align}
&\widehat{J}_1(\bus) + r_k  \widehat{J}_1'(\bus)[\bv_k] + \frac{r_k^2}{2} \widehat{J}_1''(\bus + \theta_k r_k \bv_k)[\bv_k,\bv_k]
+ \kappa g(\bus + r_kv_k)\nonumber \\
&<\widehat{J}_1(\bus) + \kappa g(\bus) +  \frac{ r_k^2}{2k}\,.\nonumber 
\end{align}
We subtract $\,\widehat{J}_1(\bus) + \kappa g(\bus)\,$ from both sides and use \eqref{directionalder} once more to find that
\begin{equation} \label{Konrad2}
r_k  \left(\widehat{J}_1'(\bus)[\bv_k]  + \kappa g'(\bus,\bv_k)\right)+ \frac{r_k^2}{2} \widehat{J}_1''(\bus + \theta_k r_k \bv_k)[\bv_k,\bv_k]<
\frac{ r_k^2}{2k}\,.
\end{equation}
From the right-hand side of \eqref{directionalder}, and the variational inequality \eqref{varineq2}, it follows that
\[
\widehat{J}_1'(\bus)[\bv_k]  + \kappa g'(\bus,\bv_k) \ge 0\,,
\]
and thus, by \eqref{Konrad2},
\begin{equation}\label{liminf1} 
\widehat{J}_1''(\bus + \theta_k r_k \bv_k)[\bv_k,\bv_k]< \frac{1}{k}\,.
\end{equation}
Passing to the limit $k \to \infty$, we  apply Lemma \ref{LA3} and deduce that \,$\widehat{J}_1''(\bus)[\bv,\bv] \le 0.$
Since we know that $\bv \in C_{\bus}$, the second-order condition \eqref{coerc} implies that \,$\bv = \mathbf{0}$. 

(iii) {\em Contradiction:} 
From the previous step we know that $\,\bv_k\to\mathbf{0}\,$ weakly in $\,L^2(Q)^2$. Moreover, \eqref{walter1} yields that
\begin{align}
\label{umba-umba}
&\widehat J_1''(\bus)[v_k,v_k] = b_3\iint_Q |\bv_k|^2\,+b_1  \iint_Q |\xi_k|^2\,+b_2\iO|\xi_k(T)|^2\nonumber\\
&+\iint_Q \Big[ 2P'(\vp^*)\xi_k(\theta_k-\chi\xi_k-\eta_k)(p^*-r^*)\,-\,F^{(3)}(\vp^*)q^*\,|\xi_k|^2 \Big]\nonumber\\
&+\iint_Q P''(\vp^*)(\s^*+\chi(1-\vp^*)-\mu^*)(p^*-r^*)|\xi_k|^2\,,
\end{align}
where we have set $\,(\eta_k,\xi_k,\theta_k)=\S'(\bus)[\bv_k]$, for $\,k\in\enne$.  
By virtue of Lemma \ref{LA3}, the sum of the last four integrals on the right-hand side converges to zero. 
On the other hand, $\,\|\bv_k\|_{L^2(Q)^2}=1\,$ for all
$\,k\in\enne$, by construction. The weak sequential semicontinuity of norms then implies that 
\begin{align*}
&\liminf_{k\to\infty} \,\widehat J_1''(\bus)[\bv_k,\bv_k] \,\ge \,\liminf_{k\to\infty} \,b_3\iint_Q|\bv_k|^2\,=\,
b_3\,>0\,.
\end{align*}
On the other hand, it is easily deduced from \eqref{liminf1} and \eqref{Lip2} that 
\[
\liminf_{k \to \infty} \widehat{J}_1''(\bus)[\bv_k,\bv_k] \le 0\,,
\]
a contradiction. The assertion of the theorem is thus proved. 
\end{proof}
\begin{remark}{\rm
We note at this place  that the formula (6.5) in \cite{CSS2}, which resembles \eqref{umba-umba},
contains three sign errors: indeed, the term in the second line of \cite[(6.5)]{CSS2} involving $ \,P''\,$ should carry a ``plus'' sign,
while the two terms in the third line should carry ``minus'' signs. These sign errors, however, do not have an impact on the validity 
of the results established in \cite{CSS2}.} 
\end{remark}


\section{Appendix}
\setcounter{equation}{0}

In the following, we assume that {\bf (A1)}--{\bf (A6)} and {\bf (C1)}--{\bf (C3)} are fulfilled and that $\bus\in\uad$ is 
fixed with associated state
$\,(\mu^*,\vp^*,\s^*)=\S(\bus)\,$ and adjoint state $\,(p^*,q^*,r^*)$. We also recall the definitions of the spaces used below 
given in \eqref{calX}, \eqref{defZ}, and \eqref{calZ}.
\Blem 
\label{LA1} 
Let $\{\bu_k\}\subset\uad$ converge strongly in $L^2(Q)^2$ to $\,\bus$, and let $\,(\mu_k,\vp_k,\s_k)=\S(\bu_k)\,$
and $\,(p_k,q_k,r_k)$, $\,k\in\enne$, denote the associated states and adjoint states. Then 
\begin{align}
\label{conmu}
&\mu_k\to\mu^*&&\quad\mbox{strongly in }\,Z,\\
\label{conphi}
&\vp_k\to\vp^*&&\quad\mbox{strongly in }\,Z\cap C^0(\overline Q),\\
\label{consigma}
&\sigma_k\to \s^*&&\quad\mbox{strongly in }\,Z,\\
\label{conp}
&	p_k\to p^*&&\quad\mbox{weakly-star in } \,Z
\,\mbox{ and strongly in }\,C^0([0,T];L^p(\Omega))\,\mbox{ for }\,1\le p <6,\\
\label{conq}
&	q_k\to q^*&&\quad\mbox{weakly-star in }\, H^1(0,T;V^*)\cap \L\infty H \cap \L2 V, 	\\
\label{conr}
&	r_k\to r^*&&\quad\mbox{weakly-star in }\,Z\, \mbox{and strongly in }\,C^0([0,T];L^p(\Omega))\,\mbox{ for }\,1\le p <6.
\end{align}
 
\Elem
\begin{proof} The strong convergence $\,\|\S(\bu_k)-\S(\bus) \|_{\cal Z}\to 0\,$ follows directly from \eqref{stabu}.
In addition, the global bound \eqref{ssbound1} implies that \,$\{\vp_k\}\,$ is bounded in the space \,$\widetilde X\,$ defined
in \eqref{calX}, which, thanks to the compactness of the embedding $\,W_0\subset C^0(\overline\Omega)$ and \cite[Sec.~8,~Cor.~4]{Simon},
is compactly embedded in $\,C^0(\overline Q)$. Therefore it holds $\,\|\vp_k-\vp^*\|_{C^0(\overline Q)}\to 0$ \,(at first only for a
suitable subsequence, but then, owing to the uniqueness of the limit point, eventually for the entire sequence). The convergence 
properties \eqref{conmu}--\eqref{consigma} of the state variables are thus shown. In addition, it immediately follows from the mean
value theorem and \eqref{ssbound3} that, as $\,k\to\infty$,
\begin{align}
\label{connon}
&\max_{i=1,2,3} \|F^{(i)}(\vp_k)-F^{(i)}(\vp^*)\|_{C^0(\overline Q)}\to 0,\nonumber\\
&\max_{i=0,1,2} \|P^{(i)}(\vp_k)-P^{(i)}(\vp^*)\|_{C^0(\overline Q)}\to 0.
\end{align}

Next, we conclude from the bounds \eqref{boundad} and \eqref{ssbound1} that there are a subsequence, which is again labeled by $\,k\in\enne$,
and some triple $\,(p,q,r)\,$ such that, as $k\to\infty$,
\begin{align}
\label{elchp}
&p_k\to p &&\quad\mbox{weakly-star in }\,Z\cap L^\infty(Q),\\
\label{elchq}
&q_k\to q&&\quad\mbox{weakly-star in }\,H^1(0,T;V^*)\cap L^\infty(0,T;H)\cap L^2(0,T;V),\\
\label{elchr}
&r_k\to r&&\quad\mbox{weakly-star in }\,Z\cap L^\infty(Q).
\end{align}
Moreover, by \cite[Sect.~8,~Cor.~4]{Simon} and the compactness of the embedding $V\subset L^p(\Omega)$ for $1\le p<6$, we
also have
\begin{equation}
p_k\to p, \quad r_k\to r, \quad\mbox{both strongly in }\,C^0([0,T];L^p(\Omega))\,\mbox{ for }\,1\le p<6.
\end{equation}    
From these estimates and \eqref{connon} we can easily conclude that, as $k\to\infty$,
\begin{align}
\label{conny}
&F''(\vp_k)q_k \to F''(\vp^*)q, \quad P(\vp_k)(p_k-r_k) \to P(\vp^*)(p-r),\nonumber\\
&P'(\vp_k)(\s_k-\chi(1-\vp_k)-\mu_k)(p_k-r_k) \to P'(\vp^*)(\s^*-\chi(1-\vp^*)-\mu^*)(p-r),
\end{align}
all weakly in $L^2(Q)$.

At this point, we consider the time-integrated version of the adjoint system \eqref{wadj1}--\eqref{wadj4} with
test functions in $L^2(0,T;V)$, written for $\,\vp_k,p_k,q_k,r_k$, $k\in\enne$. Passage to the limit as $k\to\infty$,
using the above convergence results, immediately leads to the conclusion that $\,(p,q,r)\,$ solves the time-integrated
version of \eqref{wadj1}--\eqref{wadj4} with test functions in $\,L^2(0,T;V)$, which is equivalent to saying that
$\,(p,q,r)\,$ is a solution to \eqref{wadj1}--\eqref{wadj4}. By the uniqueness of this solution, we must 
have $(p,q,r)=(p^*,q^*,r^*)$.
The convergence properties \eqref{conp}--\eqref{conr} are therefore valid for a suitable subsequence, and since
the limit is uniquely determined, also for the entire sequence. 
\end{proof}

\begin{corollary} 
\label{coroll3}
Let $\{\bu_k\}\subset\uad$ converge strongly in $L^2(Q)^2$ to $\,\bus$, and let
$\{\bv_k\}$ converge 
weakly to $\bv$ in $L^2(Q)^2$. Then 
\begin{equation}
\label{Hugo2}
\lim_{k\to\infty} \,\widehat {J}_1'(\bu_k)[\bv_k] \,=\, \widehat{J}_1'(\bus)[\bv]\,.
\end{equation}
\end{corollary}
\begin{proof} We have, with $\,\bu_k=(u_{k_1},u_{k_2})\,$ and $\,\bv_k=(v_{k_1},v_{k_2})$,
\[
\widehat{J}_1'(\bu_k)[\bv_k] = \iint_Q (-\h\,p_k + b_3 u_{k_1})v_{k_1} + \iint_Q (r_k + b_3 u_{k_2}) v_{k_2}.
\]
Owing to Lemma \ref{LA1}, we have, in particular, that  $\,\{-\h\,p_k + b_3 u_{k_1}\}\,$ and 
$\,\{r_k+b_3 u_{k_2}\}\,$ converge strongly in $\,L^2(Q)\,$ to $\,-\h\,p^*+b_3 u_1^*\,$ and $\,r^*+b_3 u_2^*$, respectively,
 whence the assertion immediately follows.
\end{proof}
\begin{lemma} \label{LA3} Let $\{\bu_k\}$ and $\{\bv_k\}$ satisfy the conditions of Corollary \ref{coroll3}, and
assume that $\,b_3=0$. Then
\begin{equation} \label{liminf0}
\lim_{k \to \infty} \widehat{J}_1''(\bu_{k})[\bv_k,\bv_k] = \widehat{J}_1''(\bus)[\bv,\bv].
\end{equation} 
\end{lemma}
\begin{proof}
Let $\bv_k=(v_{k_1},v_{k_2})$, $\bv=(v_1,v_2)$, $\,(\eta_k,\xi_k,\theta_k)=\S'(\bu_k)[\bv_k]$, and $(\eta,\xi,\theta)
=\S'(\bus)[\bv]$. Since $b_3=0$, we infer from \eqref{umba-umba} that we have to show that, as $k\to\infty$,
\begin{align}
&b_1\iint_Q |\xi_k|^2\,+\,b_2\iO |\xi_k(T)|^2\,+\iint_Q 2P'(\vp_k)\xi_k(\theta_k-\chi\xi_k-\eta_k)(p_k- r_k)\nonumber\\
&+\iint_Q \Big[P''(\vp_k)(\s_k+\chi(1-\vp_k)-\mu_k)|\xi_k|^2\,-\,F^{(3)}(\vp_k)q_k\,|\xi_k|^2 \Big]\nonumber\\
&\to \,b_1\iint_Q |\xi^*|^2\,+\,b_2\iO|\xi^*(T)|^2 \,+\iint_Q 2 P'(\vp^*)\xi^*(\theta^*-\chi\xi^*-\eta^*)(p^*-r^*)\nonumber\\
&\qquad +\iint_Q \Big[P''(\vp^*)(\s^*+\chi(1-\vp^*)-\mu^*)(p^*-r^*)|\xi^*|^2\,-\,F^{(3)}(\vp^*)q^*\,|\xi^*|^2\Big]\,,
\label{kuckuck}
\end{align}
where $(p_k,q_k,r_k)$ and $(p^*,q^*,r^*)$ are the associated adjoint states. By Lemma 4.1 and its proof, the convergence
properties \eqref{conmu}--\eqref{conr} and \eqref{connon} are valid. Moreover, we have
$$
(\eta_k,\xi_k,\theta_k)-(\eta^*,\varphi^*,\theta^*) = \left(\mathcal{S}'(\bu_k)-\mathcal{S}'(\bus)\right)[ \bv_k]\, +\,\S'(\bus)[\bv_k-\bv]\,.
$$
By virtue of \eqref{Lip1} and the boundedness of $\{\bv_k\}$ in $L^2(Q)^2$, the 
first summand on the right-hand side of this identity converges strongly to zero in $\,{\cal Z}$.
The second converges to zero weakly in $\,Z\times\widetilde X\times Z$. Hence, thanks to the compactness of the embedding
$\,Z\subset
C^0([0,T];L^p(\Omega))$\, for $1\le p<6$ (see, e.g., \cite[Sect.~8, Cor.~4]{Simon}), 
\begin{align}
\label{Hugo3}
(\eta_k,\xi_k,\theta_k)\to(\eta^*,\xi^*,\theta^*)& \quad\mbox{strongly in }\,
C^0([0,T];L^p(\Omega))^3 \quad \mbox{ for $1\le p<6$}.
\end{align}
In particular, as $k\to\infty$,
\begin{equation} \label{Hugo4}
b_1\iint_Q |\xi_k|^2\,+\,b_2\iO |\xi_k(T)|^2 \,\,\to\,\,b_1\iint_Q |\xi^*|^2\, +\,b_2 \iO |\xi^*(T)|^2\,.
\end{equation}
Moreover, owing to the strong convergences in $C^0([0,T];L^p(\Omega))$ for $1\le p<6$, it is easily checked, using H\"older's inequality, that
\begin{align}
\label{Hugo5}
&P'(\vp_k)\xi_k(\theta_k-\chi\xi_k-\eta_k)(p_k-r_k)\,\,\to\,\,P'(\vp^*)\xi^*(\theta^*-\chi\xi^*-\eta^*)(p^*-r^*)\,,\nonumber\\
&P''(\vp_k)(\s_k+\chi(1-\vp_k)-\mu_k)(p_k-r_k)|\xi_k|^2\to P''(\vp^*)(\s^*+\chi(1-\vp^*)-\mu^*)(p^*-r^*)|\xi^*|^2,
\end{align}
both strongly in $\,L^1(Q)$.
It remains to show that, as $k\to\infty$, 
$$
\iint_Q F^{(3)}(\vp_k)q_k\,|\xi_k|^2 \,\to\,\iint_Q F^{(3)}(\vp^*)q^*\,|\xi^*|^2\,.
$$
Since $\,q_k\to q^*\,$ weakly in $L^2(Q)$ by \eqref{elchq}, it thus suffices to show that \,$F^{(3)}(\vp_k)\,|\xi_k|^2 \,\to\, F^{(3)}(\vp^*)\,|\xi^*|^2\,$ 
strongly in $\,L^2(Q)$. However, this is a simple consequence of \eqref{connon} and \eqref{Hugo3}. The assertion is thus proved.
\end{proof}

\End{document}